\documentclass[dvips,ejs,preprint]{imsart}

\RequirePackage[OT1]{fontenc}
\usepackage{amsmath,amsthm,amssymb}
\usepackage[isolatin]{inputenc}
\usepackage{amsfonts,dsfont}
\usepackage{pstricks}
\usepackage{comment}
\usepackage{natbib}
\usepackage{colordvi}
\usepackage[colorlinks]{hyperref}
\usepackage{graphicx}

\doi{10.1214/09-EJS430}
\pubyear{2009}
\volume{3}
\firstpage{678}
\lastpage{711}

\startlocaldefs
\usepackage{amsthm}
\newtheorem{theorem}{Theorem}[section]
\newtheorem{proposition}[theorem]{Proposition}

\newtheorem{lemma}[theorem]{Lemma}

\newtheorem{algorithm}[theorem]{Algorithm}

\theoremstyle{definition}
\newtheorem{definition}[theorem]{Definition}

\renewcommand{\P}{\mathbb{P}}
\newcommand{\E}{\mathbb{E}} 
\newcommand{\telque}{\;|\;}
\newcommand{\FDR}{\mbox{FDR}}

\newcommand{\mbe}{\mathbb{E}}

\newcommand{\mbp}{\mathbb{P}}

\newcommand{\mtc}{\mathcal}
\newcommand{\mbf}{\mathbf}

\newcommand{\wt}[1]{{\widetilde{#1}}}
\newcommand{\wh}[1]{{\widehat{#1}}}
\newcommand{\ol}[1]{\overline{#1}}

\newcommand{\ind}[1]{{\mbf{1}\{#1\}}}
\newcommand{\e}[1]{\mbe\brac{#1}}

\newcommand{\prob}[1]{\mbp\brac{#1}}
\newcommand{\probnn}[1]{\mbp\bracnn{#1}}

\newcommand{\brac}[1]{\left[#1\right]}
\newcommand{\bracnn}[1]{\bigl[#1\bigr]}

\newcommand{\cH}{{\mtc{H}}}

\newcommand{\cS}{{\mtc{S}}}

\newcommand{\G}{{\widehat{\mathbb{G}}}}

\newcommand{\w}{\mbf{w}}
\newcommand{\W}{\mbf{W}}
\newcommand{\Pow}{\mbox{Pow}}
\newcommand{\LSU}{\mbox{\bf LSU}}

\newcommand{\SU}{\mbox{\bf SU}}
\newcommand{\SD}{\mbox{\bf SD}}
\newcommand{\fix}{\mathcal{I}}
\newcommand{\fixdown}{\mathcal{J}}
\endlocaldefs

\begin{document}

\begin{frontmatter}


\title{Optimal weighting for false discovery rate control}
\runtitle{Optimal weighting for FDR control}

\begin{aug}
\author{\fnms{Etienne} \snm{Roquain}\thanksref{t1}\ead[label=e1]{etienne.roquain@upmc.fr}}
\address{
University of Paris 6, LPMA,\\ 4, Place Jussieu, 75252 Paris cedex 05, France\\
\printead{e1}\\
\phantom{E-mail: etienne.roquain@upmc.fr\ }}
\end{aug}
\and
\begin{aug}
\author{\fnms{Mark A.} \snm{van de Wiel}
\ead[label=e2]{mark.vdwiel@vumc.nl}}
\address{
VU University \& VU University Medical Center,\\
De Boelelaan 1117, 1007MB Amsterdam, The Netherlands\\
\printead{e2}\\
}

\thankstext{t1}{Research partly carried out at the VU University, Amsterdam.}

\runauthor{E. Roquain and M. van de Wiel}
\end{aug}

\begin{abstract}
How to weigh the Benjamini-Hochberg procedure? In the context of
multiple hypothesis testing, we propose a new step-wise procedure that
controls the false discovery rate (FDR) and we prove it to be more
powerful than any weighted Benjamini-Hochberg procedure. Both
finite-sample and asymptotic results are presented. Moreover, we
illustrate good performance of our procedure in simulations and a
genomics application. This work is particularly useful in the case of
heterogeneous $p$-value distributions.
\end{abstract}



\received{\smonth{6} \syear{2009}}

\end{frontmatter}


\section{Introduction}\label{sec:intro}

In many practical situations, e.g. functional Magnetic Resonance
Imaging (fMRI) or microarray data, the problem of testing
simultaneously a large number $m$ of null hypotheses arises. Since
Neyman-Pearson's approach is the most commonly used strategy for single
testing, many researches have focused on generalizing this approach to
the multiple testing case. First, one should choose a global type I
error to be controlled, as the probability of making at least one false
discovery (family-wise error rate, FWER) or more recently the mean
proportion of false discoveries among all the discoveries (false
discovery rate, FDR, see \cite{BH1995}). Second, one should build a
procedure that controls the so-chosen global type I error rate. For
instance, Benjamini and Hochberg (\citeyear{BH1995}) proved that the linear step-up procedure (LSU)
controls the FDR when the underlying tests are independent. Third,  one
should show that the obtained procedure has good power performance, the
power being generally defined as the expected number of true
discoveries.

To our knowledge, while the two first points above are widely studied (e.g. building FDR controlling procedures,  see e.g. \cite{BY2001,Storey2002,Sar2002})), the last point is most of the time evaluated with simulations, without a full theoretical support.
 Only few works have studied rigorously the optimality of certain classes of multiple testing procedures (see \cite{LRS2005,WR2006,RDV2006,Storey2007,FDR2009}).

Maximizing the power while controlling the FDR remains a difficult task, because the FDR involves a random denominator (the number of discoveries). 
The present paper gives a contribution to the latter maximization problem, in the simple case where the null and alternative distributions are known
. This framework is natural for the power maximization of tests, as it was  also used in Neyman-Pearson's lemma for single testing. Although leading to oracle procedures, it can be used in practice as soon as the null and alternative distributions are estimated or guessed reasonably accurately from independent data.

More formally, assume that each hypothesis is tested using a test statistic, that can then be transformed into a $p$-value $p_i$, and denote by $F_i$ the alternative c.d.f. 
of $p_i$.  In general, the $F_i$'s can be possibly very different (e.g. with heterogeneous underlying data) and the $p$-values cannot be considered interchangeably. Therefore, a $p$-value weighting approach seems appropriate to improve the performance of a multiple testing procedure. This technic, that can be traced back to \cite{Holm1979}, consists in replacing in input each original $p$-value $p_i$ by the weighted $p$-value $p'_i=p_i/w_i$ for some weight vector $(w_1,\ldots,w_m)$ summing to $m$. Here, we focus on the weighted version of the LSU procedure that was proposed in \cite{GRW2006} (see also \cite{BR2008EJS}). In the latter paper, it was demonstrated that the weighted LSU still controls the FDR for any weighting (under independence between the $p$-values), and that some of these procedures could improve the power of the LSU asymptotically. In the present paper, \textit{we aim to find the most powerful procedure among all the weighted LSU procedures}, or more precisely, \textit{to find a procedure that mimics the best procedure among  the weighted LSU procedures}. Moreover, this procedure should be computable from the $p$-values distributions, i.e. the $F_i$'s.

When using the weighted version of the FWER-controlling Bonferroni procedure,  \cite{WR2006} (see also \cite{RDV2006}) have found the optimal weighting. 
 In \cite{Storey2007}, an optimal procedure was also proposed, maximizing the expected number of true discoveries while controlling the expected number of false discoveries. 
All these procedures use deterministic thresholds, which make the power maximization feasible. However, in the case of the FDR-controlling weighted LSU, the threshold depends on the final number of discoveries and a power maximization seems very difficult to make, even in the asymptotic framework where the number of $p$-values $m$ tends to infinity (see \cite{GRW2006}).

The main idea of this paper is to 
find the optimal weights simultaneously for all the possible rejection proportions $u\in [0,1]$. These multi-weights are then collected in optimal weight functions $u\mapsto W_i^{\star}(u)$ which in turn are sequentially integrated in a step-up procedure.
While the LSU procedure uses as threshold function $u\mapsto \alpha u$, we find that the new procedure uses a threshold function $u\mapsto \alpha u W_i^{\star}(u)$ not necessarily linear (and depending on the $F_i's$).

The new procedure, called ``optimal multi-weighted step-up procedure'',  will be presented in detail in Section~\ref{sec:new_proc}.
In Section~\ref{sec:main_res}, we show that it enjoys the following properties:
\begin{itemize}
\item[(i)] FDR control for a finite number of hypotheses, up to slight modifications;
\item[(ii)] power optimality for a finite number of hypotheses, up to error terms;
\item[(iii)] power optimality without error term and FDR control without modification when the number of hypotheses $m$ tends to infinity (in a specific asymptotic setting).\end{itemize}
These results are established in two different (classical) models of $p$-values, both assuming independence between the $p$-values. The results (ii) and (iii) additionally use that the $F_i$'s are strictly concave functions and that the maximization of the power at any rejection proportion is feasible, which remain quite mild assumptions.

In Section~\ref{sec:simu}, we present a simulation study which exhibits the behavior of the new procedure when the $F_i$'s are correctly specified or misspecified.  Section~\ref{sec:simu} discusses some applications and our conclusions are given in Section~\ref{sec:disc_conl}. All our results are proved in Section~\ref{technic_lemma}, while some technical parts are gathered in Appendix. Our proofs mainly use the ``self-consistency condition'' introduced in \cite{BR2008EJS} (see also \cite{FDR2009}) and Hoeffding's inequality (see \cite{Hoeff1963}).

\section{Preliminaries}\label{sec:prel}

\subsection{Models for the $p$-values}\label{sec:setting}

Let us first define the two different models for the $p$-values that will be used throughout the paper.

We consider a finite set of $m$ null hypotheses on a probability space and we let $H_i :=0$ (resp. $1$) if the $i$-th null hypothesis is true (resp. false). Letting $\mathbf{H}:=(H_i)_{ 1\leq i\leq m}\in\{0,1\}^m$, we denote by $\cH_0:=\{i\in\{1,\ldots,m\} \telque H_i=0\}$ the set corresponding to the true null hypotheses and by $m_0:=|\cH_0|$ its cardinal. Analogously, we define $\cH_1:=\{i\in\{1,\ldots,m\} \telque H_i=1\}$ and $m_1:=|\cH_1|$ for the alternative hypotheses. Since $\cH_1$ is the complement of $\cH_0$ in $\{1,\ldots,m\}$, we have $m_1=m-m_0$. The proportion of true nulls (resp. false nulls) is denoted by $\pi_0:=m_0/m$  (resp. $\pi_1:=m_1/m$) as usual.
We suppose that for the $i$-th null hypothesis it is given a \textit{$p$-value} $p_i$ i.e. a measurable function from the observation space into $[0,1]$ such that the distribution of $p_i$ is uniform on $[0,1]$ 
when the $i$-th null hypothesis is true:
\begin{equation}
\forall i\in \cH_0, \: \forall t\in[0,1], \:\P(p_i\leq t)= t.\label{equ_pval_lowbounded}\end{equation}
Under the alternative, we denote by $F_i$ the cumulative distribution function of $p_i$:  $\forall i\in \cH_1,$ $\forall t\in[0,1]$, $F_i(t):= \P(p_i\leq t)$. In our setting, the $F_i$'s are allowed to be different and we denote $\mathbf{F}:=(F_i)_{ i\in\cH_1}$ the family of alternative c.d.f.'s.  {The $p$-values are assumed mutually independent}. 
The latter model has parameters $(\mathbf{H},\mathbf{F})$ and will be referred troughout the paper as the \textit{conditional model} (because it uses a fixed vector $\mbf{H}$).

Additionally, we will consider the so-called random effects model (see e.g. \cite{ETST2001,Storey2003,GW2004}). In this model, $\mbf{H}$ is generated independently from all other random variables, from $m$ i.i.d. Bernoulli priors. The probability for a null to be true (resp. false) is denoted by $\pi_0:=\P(H_i=0) \in (0,1)$ (resp. $\pi_1:=1-\pi_0$). Then, the $p$-values are assumed to follow the conditional model conditionally to $\mbf{H}$:  the $p$-values are mutually independent conditional to $\mbf{H}$ and each $p_i$ is uniform conditional to $H_i=0$ (i.e. satisfies \eqref{equ_pval_lowbounded} conditional to $H_i=0$) and has for alternative c.d.f. $F_i$ conditional to $H_i=1$. As a consequence, unconditionally, the $p$-values are independent and for $i=1,\ldots,m$, the  c.d.f. of each $p$-value $p_i$ is
$t\mapsto\pi_0 \:t + \pi_1\: F_i(t).$
This model has for parameters $(\pi_0,\mbf{F})$ where $\mathbf{F}=(F_i)_{ 1\leq i \leq m}$ is the family of alternative c.d.f.'s. The latter model will be referred trough the paper as the \textit{unconditional
model}.\looseness=-1

\subsection{Assumptions and notation}\label{sec:notationassump}

We introduce the following possible regularity assumptions on the parameter $\mathbf{F}$ of each model, the derivative of $F_i$ being denoted by $f_i$:
\begin{align}
&\mbox{the $F_i$'s are continuous, strictly concave functions on $[0,1]$;} \tag{$A1$}\label{assump_class1}\\
&\mbox{the $F_i$'s  are twice differentiable on $(0,1)$;}\label{assump_class2}\tag{$A2$}\\
&\mbox{the functions $i\mapsto f_i(0^{+})$ and $i\mapsto f_i(1^{-})$ are constant}\label{assump_class3}\tag{$A3$}\\
&\mbox{for each $i,j$, $\lim_ {y\rightarrow f_i(0^+)}f_j^{-1}(y)/f_i^{-1}(y)$ exists in $[0,+\infty]$.}\label{assump_class4}\tag{$A4$}
\end{align}
As illustration, the assumptions \eqref{assump_class1}--\eqref{assump_class4} are all satisfied in the one-sided testing Gaussian case where we test for any $i$ the null ``$\mu_i = 0$'' against ``$\mu_i>0$'' from a Gaussian test statistic of mean $\mu_i$ and variance $1$. In that case, we have
\begin{equation}\label{equ_pvalueGauss}
 F_i(x)=\overline{\Phi}\big(\overline{\Phi}^{-1}(x)-\mu_i\big)\ \ \text{and}\ \
f_i(x)= \exp\bigg\{\mu_i\bigg(\overline{\Phi}^{-1}(x) -\frac{\mu_i}{2}\bigg)\bigg\},
\end{equation}
where we denoted $\overline{\Phi}(z):=\prob{Z\geq z}$ for $Z\sim\mathcal{N}(0,1)$.

Finally, 
for any non-decreasing function $F:[0,1]\rightarrow[0,1]$ we  denote
\begin{align*}
\fix(F)&:=\sup\{u\in[0,1] \telque F(u)\geq u\},\\
\fixdown(F)&:=\sup\{u\in[0,1] \telque \forall u'\leq u, F(u')\geq u'\},
\end{align*}
and for $\lambda>0$,
\begin{align*}
\fix_\lambda^+(F) &:=  (\fix(F)+\lambda) -  F(\fix(F)+\lambda)\:\:\: \mbox{ when } \lambda<1-\fix(F),\\
\fix_\lambda^-(F) &:= F(\fix(F)-\lambda) -  (\fix(F)-\lambda)\:\:\:\mbox{ when } \lambda<\fix(F).
\end{align*}
We easily check that $\fix(F)$ and $\fixdown(F)$ are maxima and that $F(\fix(F))=\fix(F)$ and $F(\fixdown(F))=\fixdown(F)$.

\subsection{Multiple testing procedures, FDR and power}

A \textit{multiple testing procedure} $R$ is defined as an algorithm which, from the data, aims to reject part of the null hypotheses.
Below, we will consider, as is usually the case, multiple testing procedures which can be written as a function of the family of $p$-values $\mathbf{p}=(p_i, i\in\{1,\ldots,m\})$.
More formally, we define a multiple testing procedure as a measurable function $R$, which takes as input a realization of the $p$-value family $\mbf{p}\in [0,1]^m$ and which returns a subset $R(\mbf{p})$ of $\{1,\ldots,m\}$, corresponding to the rejected hypotheses (i.e.  $i\in R(\mbf{p})$ means that the $i$-th hypothesis is rejected by $R$).

As introduced by \cite{BH1995}, the false discovery rate (FDR) of a multiple testing procedure is defined as the mean proportion of true hypotheses in the set of the rejected hypotheses:
\begin{equation}
\label{equ_FDR}
\FDR(R)=\E\bigg[\frac{|\cH_0\cap R(\mbf{p})|}{|R(\mbf{p})|\vee 1}\bigg],
\end{equation}
where $|\cdot|$ denotes the cardinality function. Of course, the FDR in \eqref{equ_FDR} depends on the model chosen for the $p$-values. In particular, the FDR in the conditional model involves an expectation taken conditionally to $\mbf{H}$, whereas the FDR in the unconditional model additionally uses an averaging over $\mbf{H}$.
It is worth noticing that, if a procedure controls the FDR in the conditional model, that is conditionally to any value of $\mbf{H}\in\{0,1\}^m$, it controls also the FDR unconditionally.

Finally, we use the standard power criterium equal to the mean proportion of correctly rejected hypotheses, that is,
\begin{equation}
\label{equ_Pow}
\Pow(R)=m^{-1} \E\big[ |\cH_1\cap R(\mbf{p})|\big].
\end{equation}

In the notation below, we will sometimes drop the explicit dependency in $\mbf{p}$ for short, writing e.g. $R$ instead of $R(\mbf{p})$.

\subsection{Weighted linear step-up procedures}\label{def_stepupdown}

Let us consider $\w=(w_i)_i$ a vector of non-negative real numbers such that $\sum_{i=1}^m w_i=m$, called here a \textit{weight vector}, and consider the weighted $p$-values $p'_{i}=p_{i}/w_i$, ordered as: $p'_{(1)} \leq \dots \leq p'_{(m)}$ with the convention $p'_{(0)}=0$.

As introduced by \cite{GRW2006},
 the \textit{weighted linear step-up procedure} associated to $\w$, denoted here by {\bf LSU}$(\w)$,  rejects the $i$-th hypothesis if $p'_i\leq \alpha \wh{u}$, with
 \begin{equation}\wh{u}=m^{-1} \max\{r\in\{0,1,\ldots,m\}\telque p'_{(r)}\leq \alpha r/m\}.\label{def_WLSU_bad}\end{equation}
In particular, the procedure {\bf LSU}$(\w)$ using $\forall i,w_i=1$ corresponds to the standard linear step-up procedure of \cite{BH1995}, denoted here by $\LSU$.
Letting $\wh{\mathbb{G}}_\w(u) = m^{-1} \sum_{i=1}^m \ind{p_i\leq \alpha w_i u}$, the rejection proportion $\wh{u}$ can equivalently be defined as
\begin{equation}\wh{u}=\fix(\wh{\mathbb{G}}_\w),\label{defWLSU} \end{equation} using the notation of Section~\ref{sec:notationassump}. %
Contrary to \eqref{def_WLSU_bad}, expression \eqref{defWLSU} does not make any specific use of the reordered $p$-values $p'_{(1)},\ldots,p'_{(m)}$, so that it is generally  more convenient from a mathematical point of view. 

For any choice of weight vector $\w$, \cite{GRW2006} proved that the weighted linear step-up procedure controls the FDR at level $\alpha m^{-1} \sum_i (1-H_i) w_i \leq \alpha$ in the conditional model and (thus) at level $\pi_0 \alpha \leq \alpha$ in the unconditional model.

\section{New approach}\label{sec:new_proc}

We present in this section a new family of multiple testing procedures, called multi-weighted procedures. 
We start by motivating their introduction from the power optimization problem among the family of weighted linear step-up procedures.

\subsection{Weight functions}\label{sec:weightchoice}

Following \cite{GRW2006}, the explicit computation of the power of the LSU$(\w)$ is a difficult task (even asymptotically): it depends on the final proportion of rejections of the procedure $\wh{u}=\fix(\wh{\mathbb{G}}_\w)$, which is a random variable 
itself depending on $\w$. Therefore, we propose here to perform the optimization for each fixed rejection proportion $u$ which in turn leads to a family of optimal weight vectors 
depending on $u$, $0< u \leq 1$. 

First, define the  power of the procedure that thresholds each $p$-value $p_i$ at level $\alpha w_i u$:
\begin{equation}\label{equ:power}
\Pow_u(\w):= \Pow(\{i\telque p_i\leq \alpha w_i u\})
,
\end{equation}
corresponding intuitively to the ``power of the LSU$(\w)$ at rejection proportion~$u$''.

Second, define a \textit{weight (vector) function}  as a function ${\W}:u\in(0,1] \mapsto \W(u)=(W_i(u))_{i} \in (\mathbb{R}^+)^m$ such that each $\W(u)$ is a weight vector, that is, $\forall u \in (0,1]$, $\sum_{i=1}^{m} W_i(u)= m$ and such that the following property holds
\begin{equation}
\forall i\in\{1,\dots,m\}, \:\:\:u\mapsto W_i(u)\: u \mbox{ is nondecreasing on $(0,1]$.}\label{equ_hypoweights}
\end{equation}
Additionally, a weight function is said \textit{continuous} if for all $i$, $u\in(0,1] \mapsto W_i(u)$ are continuous functions.

\begin{definition}
Any weight function $\W^{\star}$ solving simultaneously the maximization problems:
\begin{equation}\label{equ_optimalweightfunction}
\forall u\in (0,1], \:\: \Pow_u(\W^{\star}(u)) = \max\big\{ \Pow_u(\w), \w \mbox{ weight vector}
\big\} ,
\end{equation}
is called the \textit{optimal weight function}.
\end{definition}
Note that $\W^{\star}$ is called here abusively ``the'' optimal weight function even if it is not proved to be unique. Of course the optimal weight function depends on the model chosen for the $p$-values.
The following proposition gives (strong) sufficient conditions for existence and unicity of the optimal weight function in the different models described in Section~\ref{sec:setting}.

\begin{proposition}\label{prop:maxpow}
Assume \eqref{assump_class1}--\eqref{assump_class2}--\eqref{assump_class3} and denote the derivative of $F_i$ by $f_i$. Then the weight function $\W^{\star}$ satisfying \eqref{equ_optimalweightfunction} exists and is unique in either of the following cases:
\begin{itemize}
\item[\textbullet] In the conditional model, if $\alpha<\pi_1$, with for all $u\in (0,1],$ 
 \begin{equation}\label{equ_cor:maxpow}
W^{\star}_i(u)=  (\alpha u)^{-1} f_i^{-1}\big(y^{\star}(u)\big)\ind{H_i=1}.
\end{equation}
\item[\textbullet] In the unconditional model, with for all $u\in (0,1]$,
\begin{equation*}
W^{\star}_i(u)=  (\alpha u)^{-1} f_i^{-1}\big(y^{\star}(u)\big).
\end{equation*}
\end{itemize}
In each case, $y^{\star}(u)$ is defined as the unique element providing $\sum_{i=1}^{m}W_i^{\star}(u)=m.$
Moreover, in both models, the weight function $\W^{\star}$ is continuous, and assuming in addition \eqref{assump_class4}, the limits $W_i^{\star}(0^+)$ exist for all $i$.
\end{proposition}

The proof, which is based on similar arguments than those proposed in \cite{RDV2006} and \cite{WR2006}, is given in Section~\ref{sec:proofweightchoice}. Of course, the optimal weight function depends on the parameters of the model: on $(\mbf{H},\mbf{F})$ in the conditional model and on $\mbf{F}$ (only) in the unconditional model.

\begin{figure}[t!]
\includegraphics[scale=.3,angle=-90]{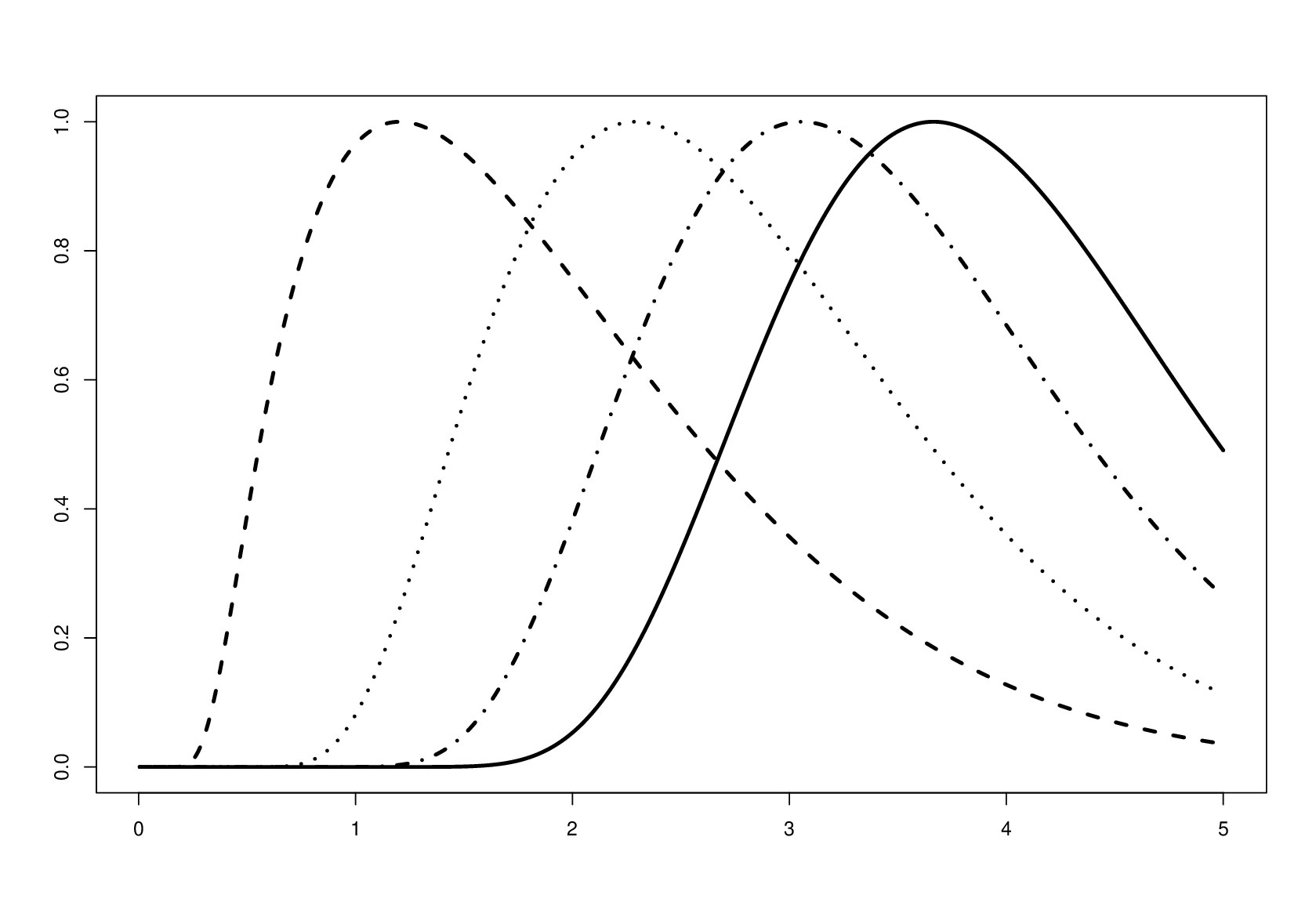}
\vspace{-3mm}
\caption{Plot of the optimal weights $(W^{\star}_i(u))_i$ in function of the alternative means $(\mu_i)_i$, for $u=1/m$ (solid), $u=10/m$ (dashed-dotted), $u=100/m$ (dotted), $u=1$ (dashed). Unconditional model and Gaussian one-sided case with $m=1000$, $\alpha=0.05$, $\mu_i=5i/m$ for $i=1,\ldots,m$. 
Each curve is normalized to have a maximum equal to $1$.}\label{fig_plotweight}
\end{figure}

For instance, when the $p$-values are generated from the Gaussian model \eqref{equ_pvalueGauss}, the optimal weight function in the unconditional model is given by
\begin{equation}\label{equ_optimal_gauss}
W_i^{\star}(u)=(\alpha u)^{-1} \ol{\Phi}\bigg(\frac{\mu_i}{2}+\frac{c(u)}{\mu_i}\bigg)
,
\end{equation}
where $c(u)$ is the unique element of $\mathbb{R}$ such that $\sum_{i=1}^{m}W_i^{\star}(u)=m.$ It therefore only depends on the vector of alternative means $\mu=(\mu_i)_{1\leq i \leq m}$.
Figure~\ref{fig_plotweight} displays the optimal weight vectors $\W(u)$ 
for a particular choice of means and different values of $u$. We observe that $\W(u)$ strongly depends on $u$:  for $u=1$, the weight vector is larger for small means, whereas as $u$ decreases, the weight vector is maximum on larger means. In particular, for small $u$, the weighting is close to zero for the smallest means, because they produce $p$-values much larger than $\alpha u$ (with high probability).

The Gaussian formula  \eqref{equ_optimal_gauss} can be also suitable for test statistics ``close to be Gaussian'', namely for locally uniform asymptotically normal test statistics (see Chapter~14 of \cite{Vaart1998} and Section~4.3 and Section~7 of \cite{RW2008}). This is the case for instance for the Mann-Whitney test statistic.

\vspace*{-3pt}
\subsection{Multi-weighted procedures}\label{section:ms}
\vspace*{-3pt}

From the previous section, we have now to integrate several weight vectors in a single multiple testing procedure, or more precisely to use a weight vector $\w$ which may depend on $u$. For this, we extend the definition of weighted linear procedures to the case of multi-weighted procedures.
\eject

First, we define the \textit{threshold collection $\Delta=(\Delta_i(u))_{i,u}$ associated to a given weight function $\mathbf{W}(\cdot)=(W_i(\cdot))_{i}$} by $\forall i\in\{1,\ldots,m\}$,
$\forall u\in [0,1] $,
\begin{center}
$\Delta_i(u):= \alpha\: W_i(u)\: u$ if $u>0$ and $\Delta_i(0):=0$.
\end{center}
Conversely, given any threshold collection $\Delta=(\Delta_i(u))_{i,u}$ such that each $\Delta_i$ is nonnegative, nondecreasing on $[0,1]$ and such that  $\forall u\in (0,1]$, $m^{-1}\sum_i\Delta_i(u)=\alpha u$, we define \textit{the weight function $\W=(W_i(u))_{i,u}$ associated to $\Delta$} by $\forall i\in\{1,\ldots,m\}$, $\forall u\in (0,1]$,
$W_i(u):=\Delta_i(u)/(\alpha u)$. As a consequence, the threshold collection $\Delta$ and the weight function $\W$ are one to one associated. 

\begin{definition} 
\label{def-mw}
Consider a weight function $\mathbf{W}(\cdot)=(W_i(\cdot))_{i}$ and its associated threshold collection $\Delta$.
 The \textit{multi-weighted step-up procedure} with weight function $\W$, denoted by \SU$(\mathbf{W})$, rejects the $i$-th null hypothesis if $p_i\leq \Delta_i(\wh{u})$, where
 \begin{equation}\label{equ_stepup_mw}\wh{u}= \fix(\wh{\mathbb{G}}_\W),\end{equation}
and where we denoted $\wh{\mathbb{G}}_\W(u) := m^{-1} \sum_{i=1}^m \ind{p_i\leq \Delta_i(u)}$ for all $u\in [0,1]$. 
 \end{definition}

In particular, in the case where for all $u$, $W_i(u)= w_i$ is independent of $u$, the procedure \SU$(\mathbf{W})$ reduces to $\LSU(\w)$. More generally, the above definition of \SU$(\mathbf{W})$ allows to choose thresholding $\Delta_i(u)$ not linear in $u$.

 As for $\LSU(\w)$, the multi-weighted procedure \SU$(\mathbf{W})$ can also be derived from a re-ordering based algorithm. The main difference is that the original $p$-values are ordered in several ways, because several weighting are used. Namely, if for $r\geq 1$, $q_r$ denotes the $r$-th smallest $\W(r/m)$-weighted $p$-value i.e. is equal to $p'_{(r)}$ where $\forall i, p'_i=p_i/W_i(r/m)$ and  putting $q_{0}=0$, we have
$$\wh{u}=m^{-1} \max\{r\in\{0,1,\ldots,m\} \telque q_r\leq \alpha r/m\}.$$

Similarly to the step-up case, we can define the \textit{multi-weighted step-down procedure} with weight function $\W$ (and associated threshold collection $\Delta$), denoted by \SD$(\mathbf{W})$, as rejecting the $i$-th null hypothesis if $p_i\leq \Delta_i(\widetilde{u})$, where
 \begin{equation}\label{equ_stepdown_mw} 
 \widetilde{u}=\fixdown(\G_\W) ,
 \end{equation}
or equivalently, $\widetilde{u}=m^{-1} \max\{r\in\{0,1,\ldots,m\} \telque $ $\forall r'\leq r, \:q_{r'}\leq \alpha r'/m\}.$

Remark that the procedures \SU$(\mathbf{W})$ and \SD$(\mathbf{W})$ only use the values of $\W(u)$ for $u\in \{1/m,2/m,\ldots,1\}$, which makes them easily computable. We refer the reader to Appendix~\ref{sec:impl} for explicit algorithmic versions of the procedures \SU$(\mathbf{W})$ and \SD$(\mathbf{W})$.

The particular multi-weighted step-up procedure \SU$(\mathbf{W}^{\star})$ using the optimal weight function $\W^{\star}$ is called the
\textit{optimal multi-weighted procedure}.
From an intuitive point of view, since this weighting maximizes the power at any rejection proportion, the latter procedure should be more powerful than any standard weighted procedure $\LSU(\w)$.
One of the goal of Section~\ref{sec:main_res} is to state this  optimality result formally.

Finally, let us remark that in the unconditional model and under the assumptions and notation of Proposition~\ref{prop:maxpow}, the optimal multi-weighted step-up procedure may be written under the following form:  reject the $i$-th hypothesis if
$f_i( p_i ) \geq y^{\star}(\wh{u}),$
where $f_i$ is the (decreasing) alternative density of $p_i$ and $y^{\star}(\wh{u})$ is adjusted from all the $p$-values, the $f_i$'s and the pre-specified level $\alpha$.
As a consequence, this procedure is based on individual tests of Neyman-Pearson's type (the observed variables being restricted to the $p$-values).


\section{Main results}\label{sec:main_res}

We present in this section the main properties of the multi-weighted procedures. First, the finite-sample FDR control of $\SU(\W)$ for any weight function $\W$, up to slight modifications. Second, the finite-sample power optimality of the procedure $\SU(\W^{\star})$ (using the optimal weight function),  up to some small error terms. Third, a consistency result, proving that the latter slight modifications are unnecessary and that error terms vanish, in a particular asymptotic setting where $m$ tends to infinity.

\subsection{Finite-sample FDR control}\label{sec:FDRcontrol}

First, let us recall  that for any choice of weight vector $\w=(w_1,\ldots,w_m)$, the weighted linear step-up procedures {\bf LSU}$(\w)$  controls the FDR at level $\alpha m^{-1} \sum_i (1-H_i) w_i \leq \alpha$ in the conditional model and at level $\pi_0 \alpha \leq \alpha$ in the unconditional model (see \cite{GRW2006,BR2008EJS}). These controls are non-asymptotic, in the sense that they are valid for any finite $m\geq 2$.

Unfortunately, the procedure \SU$(\W)$ cannot be proved to control the FDR at level $\alpha$ for any choice of weight function $\mathbf{W}$ and for any $m\geq 2$. In Appendix~\ref{appen_contrexFDRcontrol}, a (least favorable) choice of weight function is given when $m=2$, for which $\FDR(\SU(\W))$ 
slightly exceeds $\alpha$. Therefore, in order to obtain rigorous FDR control for each $m$ and any weight function, we need to slightly correct \SU$(\mathbf{W})$. 

\begin{theorem} \label{Th-FDR_indep_stepup}\label{Th-FDR_indep_stepdown}\label{Th-FDRmain}
Consider $\mathbf{W}(\cdot)=(W_i(\cdot))_{i}$ any weight function. Then for any finite $m\geq 2$, the two following procedures
\begin{itemize}
\item[\textbullet] \SU$(\widetilde{\mathbf{W}})$ with $\widetilde{W}_i(u)={W_i(u)}/({1+\alpha W_i(1)})$,
\item[\textbullet] \SD$(\widetilde{\mathbf{W}})$ with $\widetilde{W}_i(u)={W_i(u)}/({1+ \alpha u W_i(u)})$,
\end{itemize}
have their FDR less than or equal to
$$\alpha    \max_{1\leq k \leq m} \left\{m^{-1} \sum_{i=1}^m (1-H_i)  {W}_i(k/m) \right\}\leq \alpha$$
in the conditional model. As a consequence, their FDR are less than or equal to
$\alpha  \: \E \left( \max_{1\leq k \leq m} \left\{m^{-1} \sum_{i=1}^m (1-H_i)  {W}_i(k/m) \right\} \right)
\leq \alpha$ in the unconditional model.
\end{theorem}

The proof of Theorem~\ref{Th-FDR_indep_stepdown} is given in Section~\ref{proof_TH2}. Note that this result covers the earlier result of  \cite{GRW2006,BR2008EJS}, by taking $W_i(u)=w_i$ constant in $u$.

Since from \eqref{equ_hypoweights}, we have $\alpha u W_i(u)\leq \alpha W_i(1)$, both modifications of  \SU$(\mathbf{W})$ proposed above should be not too large when $\alpha W_i(1)$ is close to $0$ (e.g. when $\alpha$ is small).
Furthermore, while the correction proposed in the weighting of the step-up procedure is more conservative than the one of the step-down, a step-up procedure is always more powerful than a step-down procedure (for the same threshold collection). Therefore, in general, no modified procedure dominates  the other. 
Nevertheless, in the particular simulation setting of Section~\ref{sec:simu}, we will see that the step-down modification appears to be better.

When using the optimal weight function $\W^{\star}$, Theorem~\ref{Th-FDR_indep_stepdown} provides two modifications of the optimal multi-weighted procedure \SU$(\W^{\star})$ that control the FDR. More importantly, it shows that any misspecification in $\W^{\star}$ (e.g. in the model parameters) still leads to the correct FDR control. This is a crucial point in practice. 

Explicit finite-sample bounds for the FDR of $\SU(\W)$ -- the step-up procedure without modification -- are given in Proposition~\ref{prop_approxFDR} in the unconditional model (see Section~\ref{proof_th_asymp}). It shows that $\FDR(\SU(\W))$ should be close to $\pi_0\alpha$ when $m$ is large, so that the modifications of Theorem~\ref{Th-FDR_indep_stepup} are not needed anymore in that case. We will develop the resulting FDR consistency result more formally in Section~\ref{sec:consis} under some asymptotic conditions and for the optimal weighting.



\subsection{Finite-sample power optimality}\label{sec:pow}

For a given weight function $\W$ of associated threshold collection $\Delta$, let us denote
\begin{equation}\label{equ_moyproprej}
G_\W(u) := \e{ \widehat{\mathbb{G}}_\W(u) }= m^{-1} \sum_{i=1}^m \prob{p_i\leq \Delta_i( u)},\end{equation}
being the mean proportion of rejections at levels $(\Delta_i(u))_i$ and define similarly $G_\w(u)$ for a weight vector $\w$. Then the following theorem holds:

\begin{theorem}\label{Th-Pow-oracle}
In the unconditional model, assume that $\mbf{F}$ satisfies \eqref{assump_class1} and consider a weight function $\W^{\star}$ which maximizes the power at every proportion rejection, i.e. satisfying \eqref{equ_optimalweightfunction}.
Consider $\lambda>0$ with $\lambda<\pi_0(1-\alpha)$.
Then we have for any finite $m$,
\begin{align}
\Pow(\SU(\W^{\star})) \geq& \max_{\w} \big\{ \Pow(\LSU(\w)) -  \varepsilon (m,\fix_\lambda^+({G}_{\w}) ) \big\}  \nonumber\\
&-  \varepsilon (m,\fix_\lambda^-({G}_{\W^{\star}}) )\ind{\lambda<\fix(G_{\W^{\star}})} -   2\lambda(1- \alpha \pi_0),\label{nonasym_power_opt}
\end{align}
where $\forall x\in\mathbb{R}$,
$ \varepsilon (m,x) :=\pi_1m^2\exp\big\{-2m\big(x- m^{-1} \big)_+^2 \big\} $ and where the maximum is taken over all the weight vectors $\w$.
Moreover, we have $\fix_\lambda^-({G}_{\W^{\star}})>0$ when $\lambda<\fix(G_{\W^{\star}})$ and $\fix_\lambda^+({G}_{\w})>0$.
\end{theorem}

The proof is made in Section~\ref{proof_TH1}.
Expression \eqref{nonasym_power_opt} can be seen as a non-asymptotic ``oracle inequality'', stating that the power of the optimal multi-weighted procedure is close to
the power of the best weighted linear step-up procedure.
This finite-sample optimality result makes sense because $\SU(\W^{\star})$, as all the weighted linear step-up procedures, controls the FDR non-asymptotically at level $\alpha$ (up to the slight modifications presented in Section~\ref{sec:FDRcontrol}).

In Theorem~\ref{Th-Pow-oracle}, condition $\lambda<\pi_0(1-\alpha)$ (resp. $\lambda<\fix(G_{\W^{\star}})$)
ensures that $\fix_\lambda^+({G}_{\w}) $ (resp. $\fix_\lambda^-({G}_{\W^{\star}})$) is well defined.
Moreover, in  \eqref{nonasym_power_opt}, $\lambda$ should be chosen such that the errors terms $\varepsilon (m,\fix_\lambda^+({G}_{\w}) )$, $\varepsilon (m,\fix_\lambda^-({G}_{\W^{\star}}) )$ and $2\lambda(1- \alpha \pi_0)$  are as small as possible.
From an asymptotic point of view, assuming that the quantities $\fix_\lambda^-({G}_{\W^{\star}})$ and $\fix_\lambda^+({G}_{\w})$ are bounded away from $0$ when $m$ tends to infinity (for any fixed $\lambda$), the error terms 
tend to zero by taking successively $m$ tending to infinity and $\lambda$ tending to zero. However, the best choice $\lambda=\lambda_m$  depends
on the parameter $\mbf{F}$ and seems quite difficult to derive under an explicit form (and so are the corresponding convergence rates in  \eqref{nonasym_power_opt}).

The next section presents sufficient asymptotic conditions making $\fix_\lambda^-({G}_{\W^{\star}})$ and $\fix_\lambda^+({G}_{\w})$ bounded away from $0$ when $m$ tends to infinity, so that the error terms will asymptotically vanish in oracle inequality  \eqref{nonasym_power_opt}.



\subsection{Consistency}\label{sec:consis}

We propose in this section an asymptotic framework in which  the optimality of $\SU(\W^{\star})$ and its FDR control hold  when $m$ tends to infinity, without modification or error term.

First, we define the asymptotic setting. For all $m\geq 2$, we consider the $m$-unconditional model, where the $m$ $p$-values are chosen as the $m$ first $p$-values of an infinite sequence of independent $p$-values $(p_i)_{i\geq 1}$, each $p$-value $p_i$ having the c.d.f. $\pi_0 t + \pi_1 F_i(t)$, for a given infinite sequence of c.d.f.'s $\mathbf{F}=(F_i)_{i\geq 1}$. In this context, the weight functions depend on $m$, and we underline this dependence in the notation, by denoting $\W^{(m)}$ instead of $\W$ (and $\w^{(m)}$ instead of $\w$). 

Second, we define a \textit{converging weight function sequence} $(\W^{(m)})_m$ as a sequence of weight functions such that the associated function sequence $(G_{\W^{(m)}})_m$ (defined in \eqref{equ_moyproprej}) 
converges point-wise (on $[0,1]$).
For short, we will often use the notation $G^{\infty}$ for the limit function of $(G_{\W^{(m)}})_m$. 

\begin{theorem}\label{th_asymp}
Consider the above asymptotic framework in which $\mathbf{F}$ is assumed to satisfy \eqref{assump_class1} and consider a class $\mathcal{W}$ of converging weight function sequences.
Let $(\W^{\star, (m)})_m$ be a sequence of weight functions  such that for all $m$, $\W^{\star, (m)}$ maximizes the power at every proportion rejection  in the $m$-unconditional model (i.e. satisfies \eqref{equ_optimalweightfunction}).
For the sequence $(\W^{\star, (m)})_m$, assumed to lie in $\mathcal{W}$, and for any weight vector sequence $(\w^{(m)})_m$ belonging to $\mathcal{W}$, we additionally assume that the associated limit function  $G^{\infty}$ is continuous and satisfies
$\fix_\lambda^-(G^{\infty})>0$ for $\lambda<\fix(G^{\infty})$.
Then the multi-weighted procedure $\SU(\W^{\star,(m)})$ satisfies
\begin{eqnarray}
\lim_m\Pow(\SU(\W^{\star,(m)}))\geq  \max_{(\w^{(m)})_m}\left\{\lim_m\Pow(\LSU(\w^{(m)}))\right\},\label{asymp_optimal}
\end{eqnarray}
the maximum above being taken over any sequence of weight vectors $(\w^{(m)})_m$ belonging to $\mathcal{W}$. 
 Moreover,
 we have
 \begin{eqnarray}
\lim_m \FDR(\SU(\W^{\star,(m)})) & \leq & \pi_0 \:\alpha\,,\label{asymp_control}
\end{eqnarray}
assuming either that $\fix(G^{\infty})>0$ (with $G^{\infty}=
\lim_m G_{\W^{\star,(m)}}$) or that we have $\lim_n \lim_m m^{-1} \sum_{i=1}^m
\sup_{0<u<n^{-1}} \big\{W_i^{\star,(m)}(u)\big\}=1$.
\end{theorem}

Theorem~\ref{th_asymp} is proved in Section~\ref{proof_th_asymp}.

Under the conditions of Theorem~\ref{th_asymp}, inequalities \eqref{asymp_optimal} and \eqref{asymp_control}  imply that  $\SU(\W^{\star})$ is asymptotically more powerful than any weighted linear step-up procedure (in a certain class of converging weight vector sequences) while having the same asymptotic FDR control.
Since the uniform weighting sequence $w_i^{(m)}=1$ is always converging (with a continuous strictly concave limit function, from \eqref{assump_class1}), it can always be added in the class $\mathcal{W}$. As a consequence, the procedure $\SU(\W^{\star})$  always improves the original $\LSU$ asymptotically.
However, this should be balanced with the fact that $\SU(\W^{\star})$ uses the true parameters of the model, whereas $\LSU$ does not.

To satisfy the assumptions of Theorem~\ref{th_asymp}, we have to choose a convenient class of converging weighting sequences $\mathcal{W}$, containing the optimal weighting sequence. We give below two examples of such choice when $\mathbf{F}$  is assumed to have a particular structure.\\[-6pt]

A first example  is the case of \textit{clustered $p$-values}:
consider a parameter $\mathbf{F}$ satisfying \eqref{assump_class1} and such that $F_i$ is equal to $F_A$ (resp. $F_B$) for $i\in \cS_A^{(m)}$ (resp. $i\in\cS_B^{(m)}$), where $\{\cS_A^{(m)},\cS_B^{(m)}\}$ forms a (deterministic) partition of $\{1,\ldots,m\}$
(this model may of course be generalized to the case $K\geq 2$ clusters).
For simplicity, we  assume  that
the proportion of $p$-values $\pi_A=|\cS_A^{(m)}|/m$ in cluster $\cS_A^{(m)}$ (resp. $\pi_B=|\cS_B^{(m)}|/m$ in cluster $\cS_B^{(m)}$) does not depend on $m$ (this holds up to take a subsequence of $m$).
In this context, we merely check that a weight vector maximizing the power (at a given rejection proportion) has the same weight within a cluster. It is therefore natural to consider the following class of weighting:
\begin{align*}
\mathcal{W} = \Big\{& \bigl(\W^{(m)}\bigr)_m \:\Big|\: \forall m,  \:W_i^{(m)} = W_A  \mbox{ (resp. $W_B$) for $i\in \cS_A^{(m)}$ (resp. $i\in \cS_B^{(m)}$), }\\[3pt]
& \mbox{ for $W_A(\cdot)$,$W_B(\cdot)\geq 0$ satisfying $\forall u \in (0,1]$, $\pi_A W_A(u)+ \pi_B W_B(u) =1$, }\\[3pt]
&\mbox{ and $u\mapsto W_A(u) u$, $u\mapsto W_B(u) u$ continuous nondecreasing on $[0,1]$}\Big\}.
\end{align*}
Since for any weight function sequence $\left(\W^{(m)}\right)_m$ of $\mathcal{W}$ the function $G_{\W^{(m)}}(u)= \pi_0\alpha u + \pi_1 \pi_A F_A(\alpha u W_A(u)) +\pi_1 \pi_B F_B(\alpha u W_B(u))$ does not depend of $m$, $\mathcal{W}$ is a class of converging weight function sequences. Moreover, $G^{\infty}(u)=G_{\W^{(m)}}(u)$ is continuous and satisfies $\fix_\lambda^-(G^{\infty})>0$ for $\lambda<\fix(G^{\infty})$, either for $\W^{(m)}(u)=\w^{(m)}$ a weight vector, or for $\W^{(m)}(u)=\W^{\star,(m)}(u)$ the optimal weight function
(from one of the last statements of Theorem~\ref{Th-Pow-oracle}).
Finally, we can apply Theorem~\ref{th_asymp} to obtain the oracle inequality \eqref{asymp_optimal}. Moreover, the last assumption required for the FDR control \eqref{asymp_control}  holds assuming that the limits $W_A^{\star}(0^+)$ and $W_B^{\star}(0^+)$ exist (as is the case under assumptions \eqref{assump_class1}--\eqref{assump_class4}, see Proposition~\ref{prop:maxpow}). \\[-6pt]

A second example is the \textit{continuous one-sided Gaussian setting}, where
 for all $m$, and $1\leq i\leq m$, $F_i(x)=\overline{\Phi}\big(\overline{\Phi}^{-1}(x)-\mu(i/m)\big)$, for a mean function $\mu: [0,1]\rightarrow \mathbb{R}^+$ assumed continuous with $\mu(t)>0$ for $t>0$. In this context, we denote $F_{i/m}$ and $W_{i/m}^{(m)}$ instead of $F_{i}$ and $W_{i}^{(m)}$ for more convenience. Also note that the function $t\mapsto F_t$ can be extended to all $t$ in $[0,1]$.
 In that setting, it is relevant to consider the following class of weighting:
 \begin{align*}
\mathcal{W} = \Bigg\{& \bigl(\W^{(m)}\bigr)_m \mbox{ weight function sequence such that }\\
&\forall u\in[0,1],  \:\:\frac{1}{m} \sum_{i=1}^m F_{i/m}(\alpha u W_{i/m}^{(m)}(u))\xrightarrow[m\rightarrow \infty]{} \int_0^1F_t(\alpha u W_t(u)) dt, \\
& \mbox{ for $(W_t(\cdot))_{t\in[0,1]} \geq 0$ satisfying $\forall u$, $t\in [0,1]\mapsto W_t(u)$ continuous  }
\Bigg\}.
\end{align*}
Any weight function sequence $(\W^{(m)})_m$ of $\mathcal{W}$ is converging, because $G_{\W^{(m)}}(u)=\pi_0 \alpha u + \pi_1 {m}^{-1} \sum_{i=1}^m F_{i/m}(\alpha u W_{i/m}^{(m)}(u)) \rightarrow G^{\infty}(u) = \pi_0 \alpha u + \pi_1 \int_0^1F_t(\alpha u W_t(u)) dt$.
Moreover, expression  \eqref{equ_optimal_gauss} provides the form of the optimal weight function:
$\alpha uW_{i/m}^{\star,(m)}(u)=\ol{\Phi}\big({\mu(i/m)}/{2}+{c^{(m)}(u)}/{\mu(i/m)}\big)$ where $c^{(m)}(u)$ is taken such that $\sum_{i=1}^m W_{i/m}^{\star,(m)}(u)=m$. In Section~\ref{sec:proofweightconv}, we prove that $(\W^{\star,(m)})_m$ belongs to $\mathcal{W}$, with a ``limit weighting'' $(W^\star_t(\cdot))_{t\in[0,1]} \geq 0$ given by
$$
W_t^{\star}(u)=(\alpha u)^{-1} \ol{\Phi}\bigg(\frac{\mu(t)}{2}+ \frac{c^{\infty}(u)}{\mu(t)}\bigg),
$$
where $c^{\infty}(u)\in\mathbb{R}$ satisfies
$\int_0^1 \ol{\Phi}\big(\frac{\mu(t)}{2}+ \frac{c^{\infty}(u)}{\mu(t)}\big)dt= \alpha u$. Similarly, any weight vector sequence of the form $\w^{(m)} = \W^{\star,(m)}(u_0)$ (with $u_0$ fixed in $(0,1]$) belongs to $\mathcal{W}$  with a limit function
$G^{\infty}(u)=\pi_0\alpha u + \pi_1 \int_0^1F_t(\alpha u W_t^{\star}(u_0)) dt$ continuous and strictly concave (implying $\fix_\lambda^-(G^{\infty})>0$ for $\lambda<\fix(G^{\infty})$).
Denoting $G^{\star,\infty}$ the limit function of  $G_{\W^{\star,(m)}}$, we merely check that $G^{\star,\infty}\geq G^{\infty}$. Since this holds for any choice of $u_0$, we derive that
$\fix_\lambda^-(G^{\star,\infty})>0$ for $\lambda<\fix(G^{\star,\infty})$ (using inequalities similar to \eqref{equImoins1utile}).
 As a consequence, we may apply Theorem~\ref{th_asymp} to obtain the oracle inequality \eqref{asymp_optimal}. In particular, the power of the multi-weighted procedure $\SU(\W^{\star,(m)})$  is always asymptotically larger than the power of the weighted linear step-up $\LSU(\w^{(m)})$ for any weight vector of the form $\w^{(m)}=\W^{\star,(m)}(u_0)$, with $u_0\in(0,1]$. Roughly, the latter signifies that $\SU(\W^{\star,(m)})$ automatically chooses the best weighting among $\{\W^{\star,(m)}(u_0)\}_{u_0}$. 

\section{Simulation study}\label{sec:simu}

An important point is now to evaluate the improvement of the new multi-weighted procedure, both when we plug the true parameters or misspecified parameters in the optimal weighting. For this, we propose to perform simulations in the -- restricted but convenient -- one-sided Gaussian testing framework under the conditional model.

\subsection{Simulations framework}

We consider the problem of testing for each $i\in\{1,\ldots,m\}$, the null ``$\mu_i=0$'' against the alternative ``$\mu_i>0$'' from the observation of $m$ independent  variables $(X_i)_{i}$ with $X_i\sim \mathcal{N}(\mu_i,1)$. The parameters $(\mbf{H},\mbf{F})$ of the (conditional) model are fully determined from the vector $\mu=(\mu_i)_i$, namely by $H_i=\ind{\mu_i>0}$ and \eqref{equ_pvalueGauss}, respectively. They represent informations of a different nature:  $\mbf{H}$ provides the location of the positive means while $\mbf{F}$ supplies their values.

For all our experiments, the number of tests is $m=1000$. The vector $\mu$ is taken such that the $m_0=700$ first components of $\mu$ are equal to zero (the proportion of zeros in the mean vector is thus $\pi_0=0.7$). The $m_1=300$ remaining non-zero means are taken in two different ways:
\begin{itemize}
\item Case 1: the non-zero means increase linearly from $\frac{3}{m_1}\overline{\mu}$ to $3\overline{\mu}$.
\item Case 2: the non-zero means are gathered in three groups of different values $\overline{\mu}$, $2\overline{\mu}$ and $3\overline{\mu}$, of respective sizes $120$, $120$ and $60$.
\end{itemize}
In both cases $\overline{\mu}$ is an ``effect size'' parameter taking values in the range $0.5+0.25 k,k\in\{0,\ldots,10\}$.

The following procedures are considered:
\begin{itemize}
 \item[--]  [LSU] the linear step-up procedure $\LSU$,
  \item[--]  [LSU$^{\star}$] the step-up procedure with threshold collection $\alpha u /\pi_0$,
\item[--] [SU-W-oracle] the multi-weighted step-up procedure \SU($\widetilde{\mathbf{W}^{\star}}$) of Theorem~\ref{Th-FDR_indep_stepup}, using the optimal weight matrix $\mathbf{W}^{\star}$ (given by \eqref{equ_optimal_gauss}),
\item[--] [SD-W-oracle] the multi-weighted step-down procedure \SD($\widetilde{\mathbf{W}^{\star}}$) of Theorem~\ref{Th-FDR_indep_stepdown}, using the optimal weight matrix $\mathbf{W}^{\star}$,
 \item[--] [Unif-oracle] the weighted linear step-up procedure {\bf LSU}($\w^\star$) using a weight vector uniform on $\cH_1$: $w_i^{\star}=0$ for $\mu_i=0$ and $w_i^{\star}=m/m_1$ for $\mu_i>0$.
 \end{itemize}
The procedures [SU-W-oracle], [SD-W-oracle] and [Unif-oracle] correspond to the case where the weighting uses the true mean vector $\mu$, hence the name ``oracle''. In situations where we replace $\mu$ by a ``guess'' $\widetilde{\mu}$ in the weights, the procedures are called ``guessed'' and are denoted by [SU-W-guess], [SD-W-guess], [Unif-guess] respectively. The procedure [Unif-oracle/guess] renders a uniform weighting over the (guessed) false nulls and is close in spirit to the approach of \cite{GRW2006}. It takes only into account the subset where the hypotheses are false (``location information''), but not the values of the non-zero means.

The procedure [LSU$^{\star}$] is performed to compare with quite recent developments on $\pi_0$-adaptive procedures (see e.g. \cite{BKY2006}). Since it uses a perfect estimation of $\pi_0$, it represents the best theoretical $\pi_0$-adaptive modification of the LSU that we can build. For clarity reasons, we avoid the problem of choosing a particular estimator of $\pi_0$ and we only consider [LSU$^\star$].

All the latter procedures have provable FDR control (see Section~\ref{sec:FDRcontrol}), 
so that it is relevant to compare them in terms of power.
In all experiments the targeted FDR level is either $\alpha=0.01$ or $\alpha=0.05$.
The different performed procedures are compared in terms of relative power (RelPow) with respect to the {LSU} procedure, defined as the expected surplus proportion of correct rejections among the false nulls: for a multiple testing procedure $R$,
\begin{equation}\label{equ_relativpower}
 \mbox{RelPow}(R) := (m_1)^{-1} \big(\E (|R\cap \cH_1|)- \E (|\mbox{\bf LSU}\cap \cH_1|)\big).
\end{equation}
Roughly speaking, this relative power represents the surplus ``probability'' of a false null to be rejected with respect to
the {LSU}. This power is estimated using Monte-Carlo simulations.
Additionally, we also evaluate the ``power range'' 
defined by the power of the weighted linear procedures {\bf LSU}(${\mathbf{W}^{\star}(u_0)}$) for any $u_0\in\{1/m,2/m,\ldots,1\}$. It is represented by a gray area over the pictures.
Finally, the optimal multi-weighted step-up procedure \SU(${\mathbf{W}^{\star}}$) without correction (which controls the FDR when $m\rightarrow \infty$) is also considered, but it  is not reported on our figures, because its (relative) power is almost indistinguishable from the top of the power range.

\subsection{Procedures using the true parameters}\label{sec:simu_perfect}

\begin{figure}[t!]
\includegraphics[width=0.45\textwidth,angle=-90]{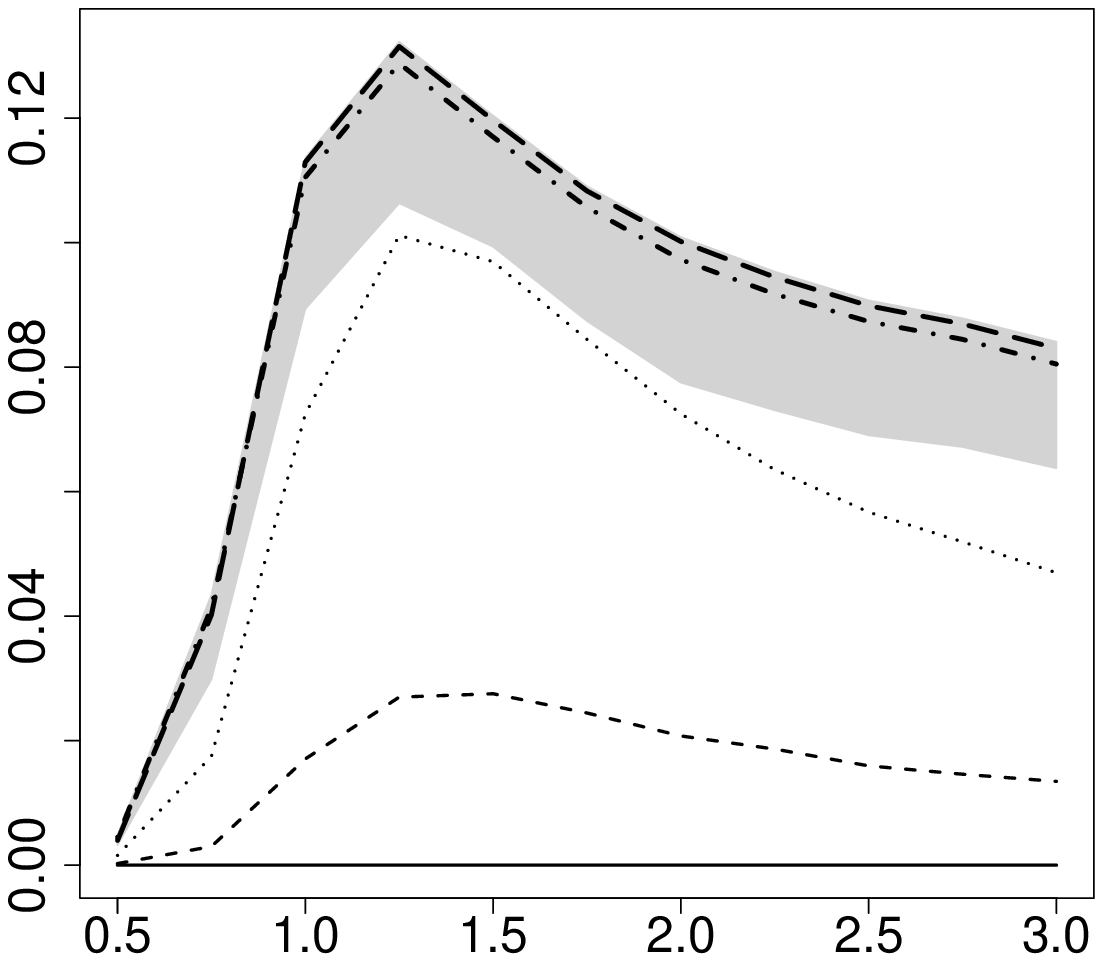} \includegraphics[width=0.45\textwidth,angle=-90]{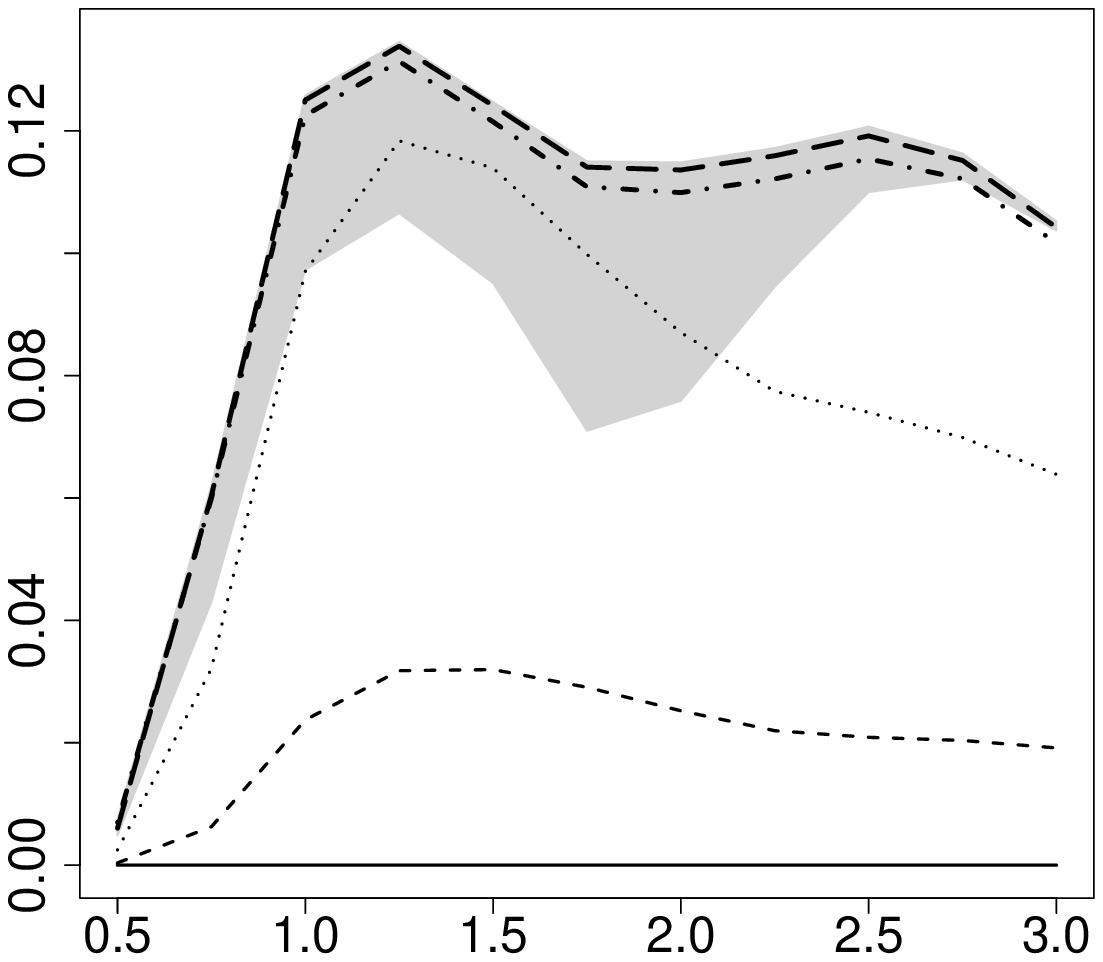}\vspace{-5mm}\\
\includegraphics[width=0.45\textwidth,angle=-90]{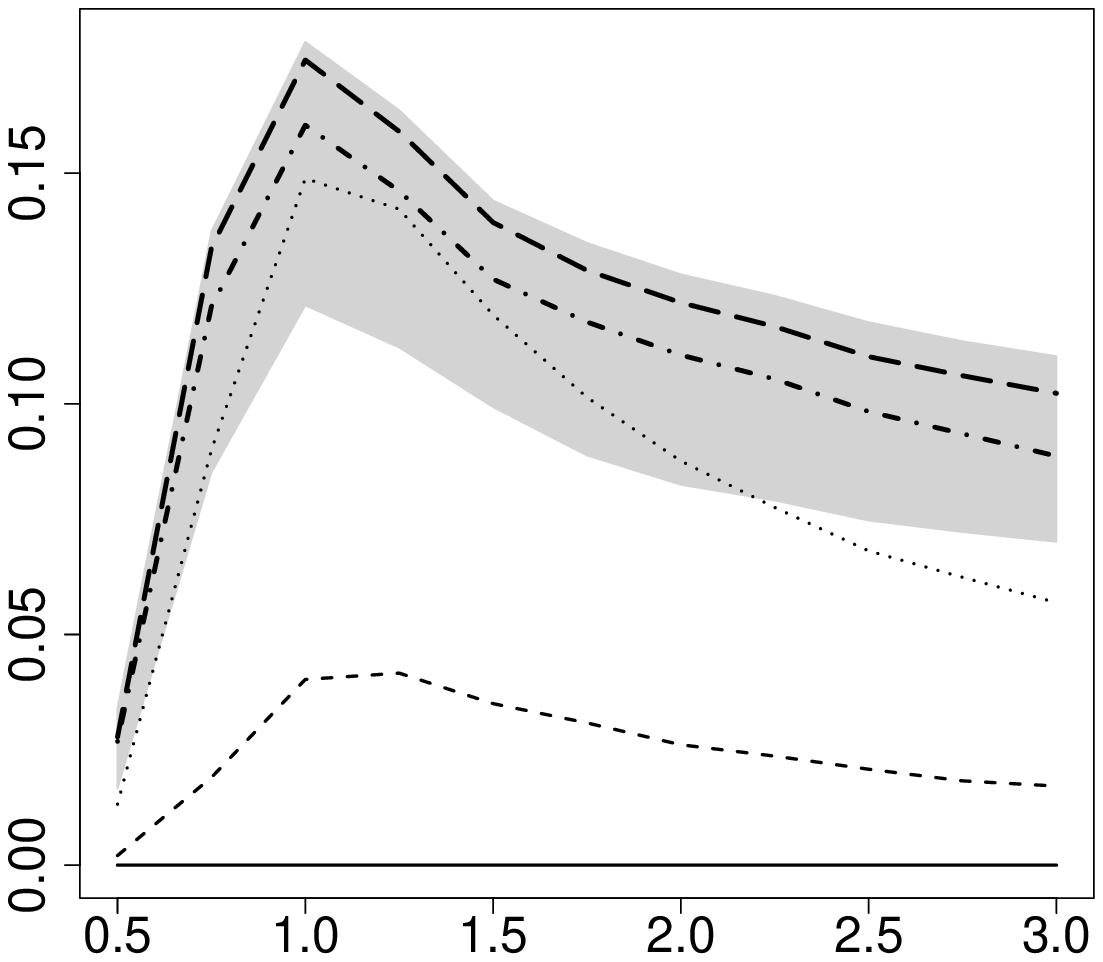} \includegraphics[width=0.45\textwidth,angle=-90]{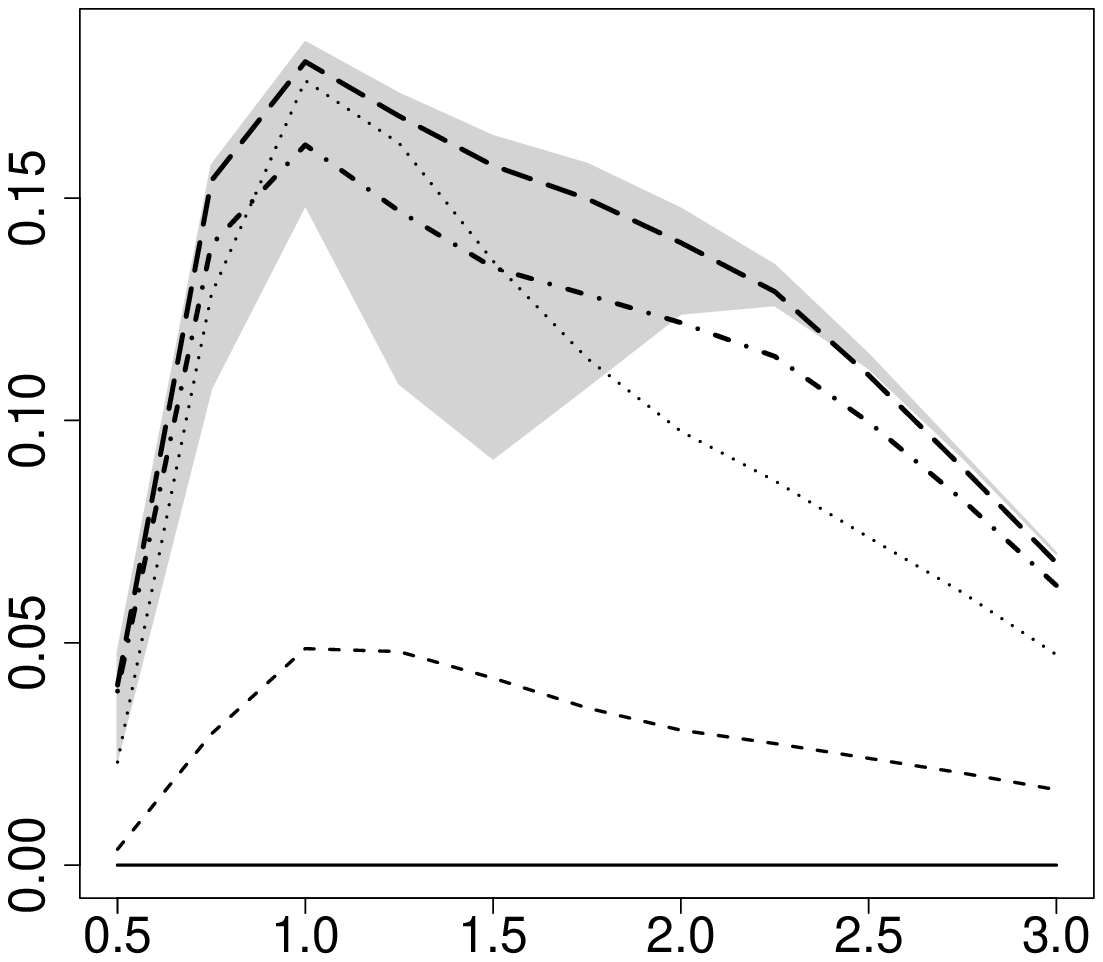}
\\
\caption{Procedures using the true parameters. Relative power of [LSU] (solid), [LSU$^{\star}$] (short-dashed), [SU-W-oracle] (dotted-dashed), [SD-W-oracle] (long-dashed) and [Unif-oracle] (dotted) in function of $\overline{\mu}$ (see text). Left: case $1$ of means, right: case $2$ of means. Top $\alpha=0.01$; bottom $\alpha=0.05$.\label{fig_powercompar_oracle}}
\vspace{-12pt}
\end{figure}

We report on Figure~\ref{fig_powercompar_oracle} the relative power \eqref{equ_relativpower} of [LSU], [LSU$^{\star}$], [SU-W-oracle], [SD-W-oracle] and [Unif-oracle] in function of the parameter $\overline{\mu}$ ($1000$ simulations). The gray area represents the power range as defined in the previous section. 

The conclusion of this experiment is that, in the most favorable case where the multi-weighting is used with the true parameters of the model, the improvement of the multi-weighted procedures over [LSU] is satisfactory. Also, [SD-W-oracle] performs here better than [SU-W-oracle] (especially for $\alpha=0.05$), so that the loss in the correction within [SU-W-oracle] seems significantly larger than the loss in the correction within [SD-W-oracle].

Furthermore, [SD-W-oracle] is more powerful than [LSU$^{\star}$] (actually, this is still true using a smaller $\pi_0$, e.g. $\pi_0=0.5$), and [SD-W-oracle] is always better than [Unif-oracle], and allows sometimes for much more discoveries. This seems coherent because [SD-W-oracle] takes into account more (correct) prior informations than [Unif-oracle]: namely, [SD-W-oracle] uses both the values and the location of the non-zero means (we are in the conditional  model), while [Unif-oracle] only uses the location information.

Finally, the procedure [SD-W-oracle] is close to the top of the power range (gray area), that is, has a power close to the power of the best procedure among {\bf LSU}(${\mathbf{W}^{\star}(u_0)}$), $u_0\in\{1/m,2/m,\ldots,1\}$. This corroborates the optimality results of Section~\ref{sec:pow}  and Section~\ref{sec:consis} in this (conditional) setting.

\subsection{Procedures using misspecified parameters}\label{sec:simu_nonperfect}

We consider here the same experiment as before, except that we take into account the ``randomness'' due to a prior guess $\widetilde{\mu}_i$ of each $\mu_i$. For this, we add a misspecification parameter $\sigma$ and we suppose that the guessed means are of the form: $\forall i\in \{1,\ldots,m\}$,
$$\widetilde{\mu}_i=\mu_i+\varepsilon_i,$$
 where $\varepsilon_i$ are i.i.d with distribution $\mathcal{N}(0,\sigma^2)$ (taken independent of the $p_i$'s). The misspecification parameter $\sigma$ is taken in the range $\{j/4,j=0,\ldots,12\}$. Remark here that the way of guessing the mean is quite raw, because it does not take into account the specific form of the parameters (of course, this guessing can be improved here by taking local means).
However, we keep this raw modeling here because we do not want to make any assumption on the parameters.

\begin{figure}[b!]
\vspace{-3mm}
\includegraphics[width=0.45\textwidth,angle=-90]{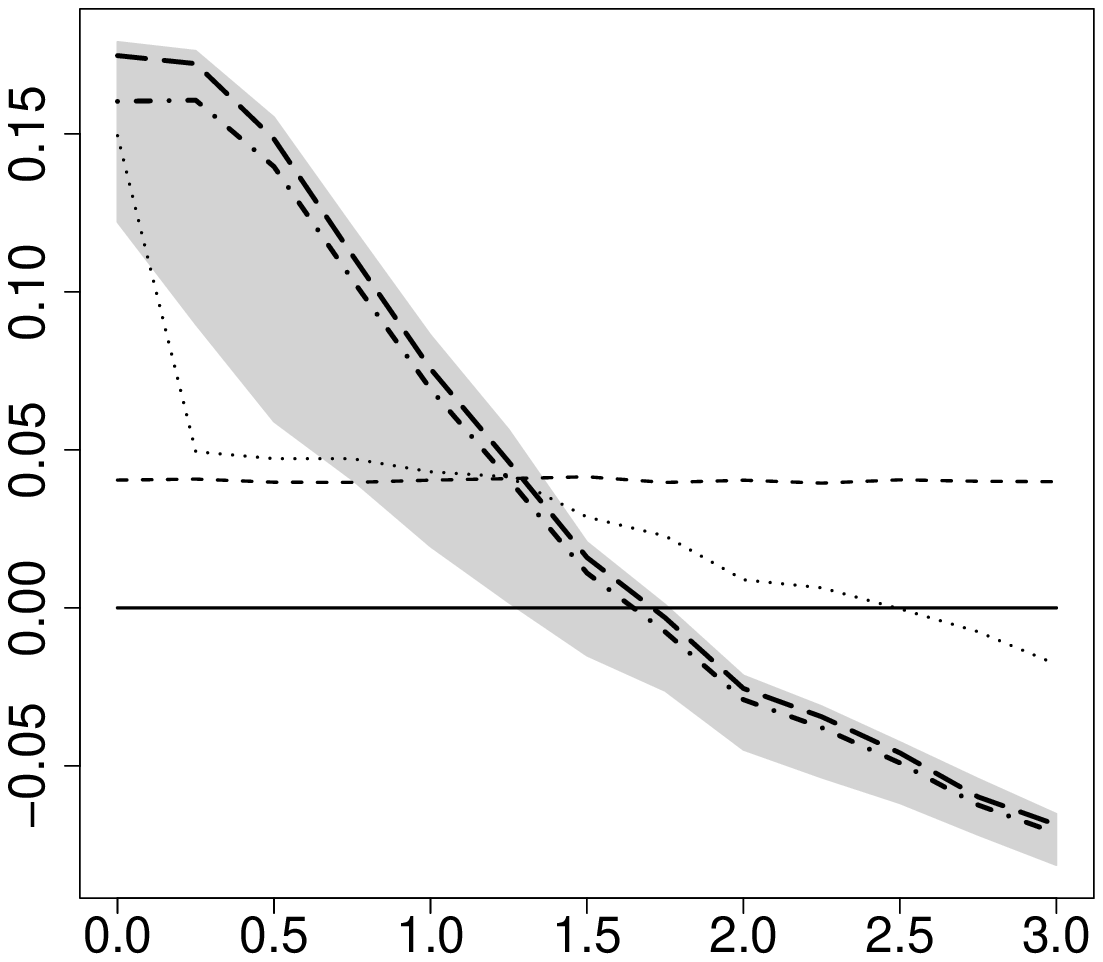}
\includegraphics[width=0.45\textwidth,angle=-90]{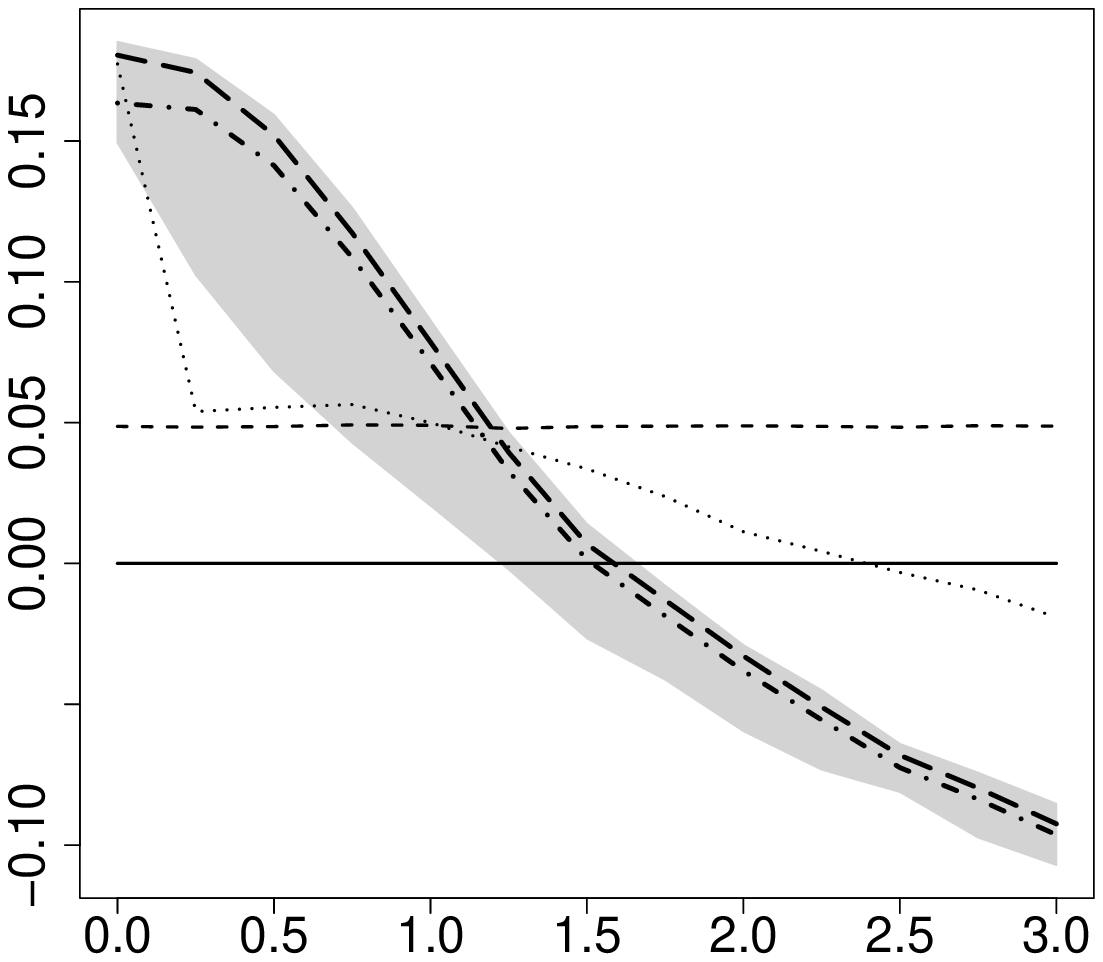}\\
\caption{Procedures using misspecified parameters. Relative power of  [LSU] (solid), [LSU$^{\star}$] (short-dashed), [SU-W-guess] (dotted-dashed), [SD-W-guess] (long-dashed) and [Unif-guess] (dotted) in function of the
misspecification parameter $\sigma$ (see text).  
Left: case $1$ of means, right: case $2$ of means (see text). $\overline{\mu}=1$;
$\alpha=0.05$.\label{fig_mispecification}}
\end{figure}

 Figure~\ref{fig_mispecification} reports the relative power \eqref{equ_relativpower} of [LSU], [LSU$^{\star}$], [SU-W-guess], [SD-W-guess] and [Unif-guess] with respect to $\sigma$.
 We performed $100$ simulations to compute the relative power and the latter is moreover averaged over $10$ generated values of the $\widetilde{\mu}_i$'s (for each values of $\sigma$).

In this experiment, we see that both multi-weighted procedures are better than other procedures when the guesses are good i.e. over the range $\sigma\in[0,1.2]$, but may be worst than the simple [LSU] procedure when $\sigma$ is large.
Furthermore, note that the procedure [Unif-guess] quickly collapses when $\sigma$ grows and therefore only proposes a slight improvement of [LSU] (or [LSU$^{\star}$]) when the guesses are good. However, it is ``less risky'' than the multi-weighted procedures for large $\sigma$.
Again, this conclusion is natural because the multi-weighted procedures take here into-account more prior information than [Unif-guess]. 

Finally, although admittedly of a limited scope, these experiments show that in principle, taking into account a correct guess of the parameters in the multi-weighted procedures should improve the power substantially. The loss/gain magnitude of these procedures depends on the quantity of prior information used (location of the positive means, value of the positive means, or both).

\section{Application to mRNA and DNA microarray experiments}\label{sec:microarray}\label{sec:appli}
In a typical microarray experiment, we want to find differentially expressed mRNA genes between two groups of individuals. For the $i$-th gene, the level expression is of the form $X_{i,1},\ldots,X_{i,k_i}$ for group 1 and $Y_{i,1},\ldots,Y_{i,\ell_i}$ for group 2, where $k_i$ (resp. $\ell_i$) is the number of individuals in group 1 (resp. group 2).

In some microarray experiments, the sample sizes $(k_i,\ell_i)$ available to assess the differential mRNA expression of gene $i$ may strongly depend on $i$, e.g. when the number of missing data differs per gene. In this application, we consider a covariate, the DNA copy number status of the same gene, which determines the groups and the sample sizes. DNA copy number status is obtained from an independent array CGH experiment, after a few pre-processing steps (see e.g. \cite{Picard2007}). We focus on the covariate $A_{i,j}$ which is equal to $1$ when gene $i$ is \textit{gained} for individual $j$ (i.e. when sample $j$ has an abnormally high DNA copy number of gene $i$), and $0$ otherwise. The biological goal behind this is to find the genes for which the mRNA expression is induced by the DNA copy number. This is particularly useful to study cancer pathologies (see e.g. \cite{Hyman2002}). Sample size dependent weights are in particular attractive here, because many genes show a large unbalance in the amount of gains (defining group 1) and non-gains (defining group 2).

Using the above framework, we analyze microarray lymphoma cancer data of \cite{Muris2007}.
In these data $m=11\,169$ genes and $n=42$ individuals are studied. The $p$-value of each gene was computed using a Mann-Whitney test. We aim to consider as prior the sample size information only, without any guess on which hypotheses are false or true. The asymptotic normality of the Mann-Whitney test statistic is used to define asymptotically optimal multi-weights $\W^{\star}$ which depend only on $(k_i,\ell_i)$ and an estimate for the global effect $\theta,$ which is a gene-independent parameter for the effect of copy number gain on mRNA gene expression. The expression of the multi-weights and the estimate for $\theta$, $\widehat{\theta}_m$, are detailed in \cite{RW2008}. The estimator $\widehat{\theta}_m$ converges in probability when $m$ grows to infinity, so that we believe that the fluctuations of $\widehat{\theta}_m$ in the weights will have a marginal effect on the effective FDR of the so multi-weighted procedure when $m$ becomes large (however we did not investigate formally the corresponding asymptotic study for now).

We applied the step-up multi-weighted procedure $\SU(\W^{\star})$, using the estimator $\widehat{\theta}_m \simeq 1.01$ of the global effect size $\theta$.
Since $m$ is large we focus on the unmodified version of our procedure, which guarantees asymptotic FDR control.
For different values of $\alpha$, the number of discoveries of this procedure and of the LSU are given in Table~\ref{tab:appli}.

\begin{table}
\caption{\label{tab:appli}Number of discoveries for the new step-up and the standard linear step-up} 
\begin{tabular}{|c||c|
c|c|}\hline
$\alpha$ & LSU & 
 $\SU(\W)$  \\\hline
 0.005 & 43 & 
   98   \\
0.01 &   121 &
 156 \\
0.05 &   478 & 
 476 \\
0.1  &  836 &
 859 \\\hline
\end{tabular}
\end{table}

We observe that our new step-up procedure discovers more differentially expressed genes when $\alpha \in \{0.005,0.01\}$. For $\alpha \in \{0.05,0.1\}$, the performance of the two step-up procedures is similar. So, the improvement of our procedure is here mostly noticeable when the proportion of rejections 
is small. This is in accordance with our intuition: the prior information (here the sample sizes) is particularly useful when the proportion of rejections is expected to be small.
Finally, let us remark that
 these positive results on the sample size problem have been corroborated in a specific simulation study as well (not reported here).

\section{Conclusion and discussions}\label{sec:disc_conl}

When the parameters of the $p$-value model are known, we proposed to solve the problem of the LSU optimal weighting by finding a new procedure which provably outperforms all the weighted LSU procedures (up to small error terms) and which can be easily computed from these parameters.
Our simulations illustrated the strength of the improvement of our new approach in situations where it uses the true or misspecified  parameters.

In our results, the assumptions concerning the marginal distributions of the $p$-values were quite mild: the FDR control only required that each $p$-value is uniform under the null while the optimality results only required the strict concavity of the c.d.f.'s of the $p$-values. Moreover, the existence of the optimal weight function only asked to maximize simultaneously the power at any proportion rejection, and we gave strong sufficient assumptions for its existence and unicity.

Several extensions to this work are possible: first, we have supposed the independence between the $p$-values  all along the paper, which is a standard but somewhat unrealistic assumption for the applications. In \cite{RW2008}, we proposed some extensions of the present FDR control results to the case of positively regressively dependence or unspecified dependence. However, the so-derived procedures seemed too conservative for practical use. Therefore, there is a room left for future investigations, which join the very active (but  challenging) research field studying the impact of $p$-value dependence on FDR control (see e.g. \cite{KW2008,RSW2008}).

Second, our FDR controls are done at level smaller than $\pi_0\alpha$ (asymptotically, in the unconditional model). Therefore, when $\pi_0$ is small, our procedures are inevitably conservative, because their actual FDR is much lower than the fixed target. This is a classical problem for the LSU procedure and several works have been proposed to address this issue, by integrating a $\pi_0$-estimate in the threshold, building so-called \textit{adaptive} LSU procedures (see e.g. \cite{BKY2006,BR2008b}). A possible interesting extension of our work could therefore be to derive adaptive multi-weighted procedures, which would increase the power when the data contain a lot of signal.

A third -- and maybe more important -- direction for future works is the investigation of data-driven weighting. A first idea could be to replace the function $\Pow_u(\cdot)$, the power at rejection proportion $u$, by an empirical substitute and to perform the simultaneous maximization with this substitute. This would yield an empirical optimal weight function $\wh{\W^\star}$ that can in turn be integrated in a multi-weighted procedure. While this certainly requires to use a model with some replications, the theoretical FDR control and power optimality of such data-driven procedure are not straightforward from the present work, because all our proofs here use the fact that the weight functions are deterministic.

\section{Proofs}\label{technic_lemma}

\subsection{Useful notation and lemmas}\label{notproof}

Let us first introduce the following notation that will be useful throughout our proofs: if $R$ is the step-up procedure associated to a given weight function $\W$ of associated threshold collection
$\Delta$, and $\wh{u}:= |R|/m$ its rejection proportion, that is $\wh{u} = \fix(\G_{\W})$, we denote by:
\begin{enumerate}
\item $R_{-i}$ the step-up procedure on the set of hypotheses corresponding to $\{1,\ldots,m\}\backslash \{i\}$, that is excluding the $i$-th null,
and associated to the threshold collection $\forall j\neq i$, $\forall u$, 
$\Delta_j( (1-m^{-1})u);$ and we denote by $\wh{u}_{-i} :=|R_{-i}|/(m-1)$ its rejection proportion, so that  $\wh{u}_{i}=\fix( \wh{\mathbb{G}}_{-i} )$ with  $\wh{\mathbb{G}}_{-i}(u) :=(m-1)^{-1}   \sum_{j\neq i} \ind{p_j\leq \Delta_j( (1-m^{-1})u)}$;
\item $R'_{-i}$ the step-up procedure on the set of hypotheses excluding the $i$-th null associated to the threshold collection $\forall j\neq i$,  $\forall u$, 
$\Delta_j((1-m^{-1})u+m^{-1});$  and we denote $\wh{u}'_{-i} :=|R'_{-i}|/(m-1)$ its rejection proportion, hence  $\wh{u}'_{-i}=\fix( \wh{\mathbb{G}}'_{-i} )$ with  $\wh{\mathbb{G}}'_{-i}(u) :=(m-1)^{-1}   \sum_{j\neq i} \ind{p_j\leq \Delta_j((1-m^{-1})u+m^{-1})}$.
\end{enumerate}
Similarly, when $R$ is step-down, we define $R_{-i}$ and $R'_{-i}$ as step-down procedures and we denote $\wt{u}:= \fixdown(\G_{\W})$, $\wt{u}_{-i}:= \fixdown( \wh{\mathbb{G}}_{-i} )$, $\wt{u}'_{-i} :=\fixdown( \wh{\mathbb{G}}'_{-i} )$ instead of $\wh{u}$, $\wh{u}_{-i} $, $\wh{u}'_{-i} $,  respectively.

The two following lemmas make a link between the rejection proportions of  $R$, $R_{-i}$ and $R'_{-i}$, for different values of $p_i$.
They are proved in Appendix~\ref{sec:lemma_proof} and are related to Lemma~10.20 of \cite{Roq2007}.

\begin{lemma} \label{lemma_stepupdown}
Let $R$ be the step-up procedure associated to a given weight function of threshold collection
$\Delta$ and consider $\wh{u}$, $\wh{u}_{-i}$ and $\wh{u}'_{-i}$ as above. Then we have point-wise:
\begin{enumerate}
\item $p_i \leq \Delta_i(\wh{u}) \:\: \Longleftrightarrow\:\: p_i \leq \Delta_i((1-m^{-1})\wh{u}'_{-i}+m^{-1})  \:\: \Longleftrightarrow\:\:  \wh{u} = (1-m^{-1})\wh{u}'_{-i}+m^{-1} $\,; 
\item $p_i > \Delta_i(\wh{u})  \:\:\Longleftrightarrow\:\:  \wh{u} =   (1-m^{-1})\wh{u}_{-i}$\,. 
\end{enumerate}
\end{lemma}

\begin{lemma} \label{lemma_stepdown}
Let $R$ be a step-down procedure associated to a given weight function of threshold collection
$\Delta$ and consider $\wt{u}$ and $\wt{u}_{-i}$ as above.  Then we have point-wise, for any $k\in\{1,\ldots,m\}$,
\begin{enumerate}
\item $\wt{u}\geq k/m$ and $p_i>\Delta_i((k-1)/m)$  $\:\:\Longrightarrow\:\:$ $\wt{u}_{-i} \geq (k-1)/(m-1)$\,;
\item $\wt{u}_{-i} \geq (k-1)/(m-1)$ and $p_i\leq \Delta_i((k-1)/m)$ $\:\:\Longrightarrow\:\:$ $\wt{u}\geq k/m$\,;
\item $p_i>\Delta_i((1-m^{-1})\wt{u}_{-i}+m^{-1} )$ $\:\:\Longrightarrow\:\:$ $ \wt{u} =   (1-m^{-1})\wt{u}_{-i}$\,.
\end{enumerate}
\end{lemma}

\subsection[Proof of Theorem~4.1 -- step-up part]{Proof of Theorem~\ref{Th-FDR_indep_stepup} -- step-up part}\label{proof_TH1}

The inequalities are established in the conditional model (the result in the unconditional model directly follows).

We use in all our FDR bounds that a procedure $R$ satisfying
the ``self-consistency condition'' $R=\{i\telque p_i\leq \Delta_i(|R|/m)\}$ has a FDR equal to
\begin{eqnarray}
\FDR(R)= \E\bigg[\frac{|R\cap\cH_0|}{|R|\vee 1}\bigg]= \sum_{i=1}^m (1-H_i)\E\bigg[\frac{\ind{p_i\leq \Delta_i(|R|/m)}}{|R|}\bigg]
 \label{equ_FDR_proof}.
\end{eqnarray}

Now, consider the multi-weighted step-up procedure $R=\SU(\wt{\W})$ of Theorem~\ref{Th-FDR_indep_stepup}, and denote by $\Delta$ the threshold collection associated to $\wt{\W}$: $\Delta_i(k/m)= \alpha \widetilde{W}_i(k/m) k/m =  \alpha {W}_i(k/m) k/m (1+\alpha W_i(1))^{-1} \leq 1$. Since any step-up procedure satisfies the self consistency condition, we may use \eqref{equ_FDR_proof}. Furthermore, using the notation of Section~\ref{notproof} and applying Lemma~\ref{lemma_stepupdown} (first statement), the assertion $p_i\leq \Delta_i(|R|/m)=\Delta_i(\wh{u}) $ is equivalent to $\wh{u} = (1-m^{-1})\wh{u}'_{-i}+m^{-1}$. Thus, we may rewrite the FDR as follows:
\begin{align}
\FDR(R) &=\sum_{i=1}^m (1-H_i)\E\bigg[\frac{\ind{p_i\leq \Delta_i(\wh{u})}}{\wh{u}\:m}\bigg] \nonumber\\
 &=\sum_{i=1}^m (1-H_i) \sum_{k=1}^{m} {k}^{-1} \prob{p_i\leq \Delta_i(k/m)  , \wh{u} \:m=k}\nonumber\\
  &=\sum_{i=1}^m (1-H_i) \sum_{k=1}^{m} {k}^{-1} \prob{p_i\leq\Delta_i(k/m)  ,(m-1)\wh{u}'_{-i}+1=k}\nonumber.
  \end{align}
  Then, since $\wh{u}'_{-i}$ only depends on the $p$-values of $(p_j,j\neq i)$, it is independent of  $p_i$ and we obtain
\begin{align}
\label{equ_boundFDRMWLSU}
\FDR(R)  &= \frac{\alpha}{m}\sum_{i=1}^m (1-H_i) \sum_{k=1}^{m} \widetilde{W}_i(k/m) \prob{(m-1)\wh{u}'_{-i}+1=k}\nonumber\\
&=  \frac{\alpha}{m}\sum_{i=1}^m (1-H_i) \sum_{k=1}^{m} {W}_i(k/m)(1-\Delta_i(1)) \prob{(m-1)\wh{u}'_{-i}+1=k}\nonumber\\
&=  \frac{\alpha}{m}\sum_{i=1}^m (1-H_i) \sum_{k=1}^{m} {W}_i(k/m) \prob{p_i>\Delta_i(1) ,(m-1)\wh{u}'_{-i}+1=k},
\end{align}
where we used that $p_i$ has a uniform distribution on $[0,1]$ (from \eqref{equ_pval_lowbounded}).
Next, consider the threshold collection $\forall j\in\{1,\ldots,m\}$, $\forall u$,  $\Delta'_j(u)= \Delta_j((1-m^{-1})u+m^{-1})$ and the associated step-up procedure that we denote by $R'$. Let us also denote its rejection proportion by $\wh{u}'=|R'|/m$.
From the definition of Section~\ref{notproof}, the restriction of $R'$ to the hypothesis set corresponding to $\{1,\ldots,m\}\backslash \{i\}$ is exactly the procedure $R'_{-i}$.
Therefore, from Lemma~\ref{lemma_stepupdown} (second statement applied to $R'$), the condition $p_i>\Delta_i(1)=\Delta'_i(1)\geq \Delta'_i(\wh{u}')$ implies $m\wh{u}' = (m-1)\wh{u}'_{-i}$. Therefore,
\begin{align}
\FDR(R) &\leq \frac{\alpha}{m}\sum_{i=1}^m (1-H_i) \sum_{k=1}^{m} {W}_i(k/m) \prob{p_i>\Delta_i(1) , m\wh{u}'+1=k }\nonumber\\
 &\leq \frac{\alpha}{m} \sum_{k=1}^{m} \Bigg[ \sum_{i=1}^m (1-H_i)  {W}_i(k/m)\Bigg]  \prob{m\wh{u}'+1=k}\nonumber\\
&\leq \alpha    \max_{1\leq k \leq m} \left\{m^{-1} \sum_{i=1}^m (1-H_i)  {W}_i(k/m) \right\}.\nonumber
\end{align}

\allowdisplaybreaks

\subsection[Proof of Theorem~4.1 -- step-down part]{Proof of Theorem~\ref{Th-FDR_indep_stepdown} -- step-down part}\label{proof_TH2}

Again, it is sufficient to look at the conditional model. First, let us prove that for any step-down procedure $R$ with threshold collection $\Delta$ and rejection proportion $\wt{u}$, we have for any $i$,
\begin{equation}\label{equ_GBS}
\FDR(R)\leq \sum_{i=1}^m (1-H_i)\sum_{k=1}^m \frac{1}{k} \prob{(m-1)\wt{u}_{-i}=k-1,p_i\leq \Delta_i( k/m)},
\end{equation}
where $\wt{u}_{-i}$ is the rejection proportion of the step-down procedure associated to $\Delta$ and restricted to the hypotheses different from the $i$-th hypothesis as defined in Section~\ref{notproof}.
This result has been implicitly proved in \cite{GBS2008} (Section~3), using a specific non-weighted step-down procedure.
Here, we state \eqref{equ_GBS} in a more general framework. Applying the two first points of Lemma~\ref{lemma_stepdown}, we obtain the following relations: 
\begin{align*}
\sum_{k=1}^m \frac{1}{k} &\prob{m\wt{u}=k,p_i\leq \Delta_i(k/m)}\\
=&\:\:\sum_{k=1}^m\frac{1}{k}\Big[\prob{m\wt{u}=k,p_i\leq \Delta_i((k-1)/m)} \\
&+ \prob{m\wt{u} = k,\Delta_i((k-1)/m)< p_i\leq \Delta_i(k/m)}\Big] \\
=&\:\: \sum_{k=1}^m \frac{1}{k} \prob{m\wt{u}\geq k,\Delta_i((k-1)/m)< p_i\leq \Delta_i(k/m)}  \\
&-\:\:\sum_{k=1}^m \bigg[\frac{\ind{k>1}}{k-1}-\frac{1}{k}\bigg] \prob{m\wt{u} \geq k, p_i\leq \Delta_i((k-1)/m)}\\
\leq&\:\: \sum_{k=1}^m \frac{1}{k} \prob{(m-1)\wt{u}_{-i}\geq k-1,\Delta_i((k-1)/m)< p_i\leq \Delta_i(k/m)} \\
&-\:\: \sum_{k=1}^m \bigg[\frac{\ind{k>1}}{k-1}-\frac{1}{k}\bigg] \prob{(m-1)\wt{u}_{-i}\geq k-1, p_i\leq \Delta_i((k-1)/m)}\\
=& \:\:\sum_{k=1}^m \frac{1}{k} 
\prob{(m-1)\wt{u}_{-i}=k-1,p_i\leq \Delta_i(k/m)}.
\end{align*}
As a consequence, the latter combined with \eqref{equ_FDR_proof} states \eqref{equ_GBS}. 

Now, consider the step-down procedure $R$ of Theorem~\ref{Th-FDR_indep_stepdown}, that is, associated to the threshold collection $\Delta_i(k/m)= \alpha \widetilde{W}_i(k/m) k/m =  \alpha {W}_i(k/m) k/m (1+\alpha W_i(k/m)k/m)^{-1} \leq 1$.
We use the independence between the $p$-values and \eqref{equ_GBS} to show
\begin{align}
\FDR(R)&\leq \sum_{i=1}^m (1-H_i)\sum_{k=1}^m \frac{1}{k} \prob{(m-1)\wt{u}_{-i}=k-1,p_i\leq \Delta_i(k/m)}\nonumber\\[2pt]
&= \frac{\alpha}{m} \sum_{i=1}^m (1-H_i)\sum_{k=1}^m W_i(k/m)(1-\Delta_i(k/m)) \prob{(m-1)\wt{u}_{-i}=k-1}\nonumber\\[2pt]
&= \frac{\alpha}{m} \sum_{i=1}^m (1-H_i)\sum_{k=1}^m W_i(k/m) \prob{(m-1)\wt{u}_{-i}=k-1,p_i>\Delta_i(k/m)}\nonumber.
\end{align}
The third point of Lemma~\ref{lemma_stepdown} thus implies
\begin{align*}
\FDR(R)&\leq \frac{\alpha}{m} \sum_{i=1}^m (1-H_i)\sum_{k=1}^m W_i(k/m) \prob{m\wt{u}=k-1,p_i>\Delta_i(k/m)}\\[2pt]
&\leq \frac{\alpha}{m} \sum_{k=1}^m \bigg[\sum_{i=1}^m (1-H_i)W_i(k/m)\bigg]\prob{m\wt{u}=k-1}\\[2pt]
&\leq \alpha    \max_{1\leq k \leq m} \left\{m^{-1} \sum_{i=1}^m (1-H_i)  {W}_i(k/m) \right\}.
\end{align*}

\subsection[Proof of Theorem~4.2]{Proof of Theorem~\ref{Th-Pow-oracle}}\label{proof_TH3}

 Let us assume that the following proposition holds (the proof is given at the end of this section):

\begin{proposition}\label{prop_approxpower}
In the unconditional model, consider a weight function $\W$ with its associated threshold collection $\Delta$ and put $\bar{u}=\fix(G_\W)$. Then the following holds:
\begin{itemize}
\item[(i)] assuming that for all $u' > u > \bar{u}$, $u'-G_\W(u')>  {u} - G_\W( {u})$, we have for all $\lambda>0$,
$\lambda<1-\bar{u}$,\vadjust{\eject}
\begin{eqnarray} \label{approxpower1}
\hspace*{-18pt}&&\Pow(\SU(\W)) -(1-\alpha \pi_0)\bar{u} \nonumber\\[-2pt]
\hspace*{-18pt}&&\leq \pi_1m^2\exp\big\{-2m\big(\fix_\lambda^+(G_\W)- m^{-1} \big)_+^2 \big\}- \fix_\lambda^+(G_\W) +  \lambda(1- \alpha
\pi_0);\quad \quad
\end{eqnarray}
\item[(ii)] assuming $\Delta\leq 1$, we have for all $\lambda>0$, $\lambda<\bar{u}$,
\begin{eqnarray}\label{approxpower2}
\hspace*{-18pt}&&\Pow(\SU(\W)) - (1-\alpha \pi_0)\bar{u}\nonumber\\[-2pt]
\hspace*{-18pt}&&\geq -\pi_1m\exp\big\{-2m(\fix_\lambda^-({G}_{\W}) )_+^2 \big\}  + \fix_\lambda^-({G}_{\W})  -   \lambda(1- \alpha \pi_0).
\end{eqnarray}
\end{itemize}
\end{proposition}

We now prove Theorem~\ref{Th-Pow-oracle} by applying Proposition~\ref{prop_approxpower}.
First, remark that, in the unconditional model, we have for any weight vector $\w$,
$$
G_\w(u) =  \alpha \pi_0 u + \Pow_u(\w),
$$
so that maximizing in $\w$ the power at rejection level $u$ is equivalent to maximize $G_\w(u)$ in $\w$. As a consequence, taking the optimal weight function $\W^\star$, we deduce from \eqref{equ_optimalweightfunction} that for any weight vector $\w$ and for any $u$ we have $G_\w(u)\leq G_{\W^{\star}}(u)$. Denoting $u_\w:=\fix(G_\w)$  and $u^{\star}:=\fix(G_{\W^{\star}})$, this in turn implies that
\begin{equation}
u_\w\leq u^{\star}
\label{equ_puissancelim_maj}
\end{equation}
Second, remark that $\W^{\star}$ has a threshold collection $\Delta^{\star}$ satisfying $\Delta^{\star}\leq 1$. The latter holds because the $F_i$'s are increasing (as non-decreasing strictly concave functions), and because
$\Delta^{\star}_i(u)\leq \alpha W^{\star}_i(1)$ with $\W^{\star}(1)$ maximizing the power at rejection proportion $1$.
Third, we check the assumption of (i) Proposition~\ref{prop_approxpower} for $\W(\cdot)$ constantly equal to  a weight vector $\w$, which directly follows from the strict concavity of $G_\w$ (itself coming from the strict concavity of the $F_i$'s).
Forth, let us prove that $\fix_\lambda^+({G}_{\w})>0$ and $\fix_\lambda^-({G}_{\W^{\star}})>0$. The first statement comes from the definition of $\fix(G_\w)$. To prove the second statement, consider $u^{\star}=\fix(G_{\W^{\star}})$ and the weight vector  $\w=\W(u^{\star})$, so that $u^{\star}$ is equal to $u_\w=\fix(G_\w)$ (because $u_\w \leq u^{\star}$ from \eqref{equ_puissancelim_maj}). Using again that $\W^{\star}$ is a maximum,
we obtain
\begin{align}
G_{\W^{\star}}(u^{\star}-\lambda ) &\geq G_{\w}(u^{\star}-\lambda )
= G_{\w}(u_\w-\lambda ) = \frac{G_{\w}(u_\w-\lambda )}{ u_\w-\lambda } (u_\w-\lambda) \nonumber\\[-2pt]
&> \frac{G_{\w}(u_\w)}{ u_\w } (u_\w-\lambda ) = u^{\star}-\lambda\,,\label{equImoins1utile}
\end{align}
by the strict concavity of $G_{\w}$. This implies $\fix_\lambda^-({G}_{\W^{\star}})>0$.

Finally,  using  \eqref{approxpower1} with $\W(\cdot)$ constantly equal to  any weight vector $\w$, together with
  \eqref{approxpower2} used with $\W=\W^{\star}$, we obtain for all $\lambda>0$, $\lambda<u^{\star}$ and $\lambda<\pi_0(1-\alpha)$,
\begin{align*}
&\Pow(\SU(\W^{\star}))\\ &\geq (1-\alpha \pi_0)u^{\star} - \varepsilon (m,\fix_\lambda^-({G}_{\W^{\star}}) ) -   \lambda(1- \alpha \pi_0)  \\
&\geq (1-\alpha \pi_0)u_\w -  \varepsilon (m,\fix_\lambda^-({G}_{\W^{\star}}) ) -   \lambda(1- \alpha \pi_0)\\
&\geq  \Pow(\LSU(\w)) -  \varepsilon (m,\fix_\lambda^+({G}_{\w}) ) -  \varepsilon (m,\fix_\lambda^-({G}_{\W^{\star}}) ) -   2\lambda(1- \alpha \pi_0),
\end{align*}
which proves \eqref{nonasym_power_opt}.

Let us now prove Proposition~\ref{prop_approxpower}. 
Using that the procedure $R=\SU(\W)$ satisfies the self-consistency condition $R=\{i\telque p_i\leq \Delta_i(|R|/m)\}$, Lemma~\ref{lemma_stepupdown} (first statement) and the notation of Section~\ref{notproof}, the power of  the procedure $R$ may be expressed as follows:
\begin{align}
\Pow(\SU(\W)) &=  \pi_1 m^{-1}\sum_{i=1}^m  \prob{p_i\leq \Delta_i(|R|/m) \telque H_i = 1} \nonumber\\[2pt]
 &= \pi_1 m^{-1}\sum_{i=1}^m   \prob{p_i\leq\Delta_i((1-m^{-1})\wh{u}'_{-i}+m^{-1})\telque H_i = 1}\nonumber\\[2pt]
 &\leq \pi_1 m^{-1}\sum_{i=1}^m \e{F_i \circ \Delta_i((1-m^{-1})\wh{u}'_{-i}+m^{-1})}\label{equ_interm_inequ3},
\end{align}
where we used both the independence between $p_i$ and $\wh{u}'_{-i}$ conditionally to $\mbf{H}$ and the independence between $\wh{u}'_{-i}$ and $H_i$.
For simplicity, we introduce the increasing function $\phi(u):=(1-m^{-1})u+m^{-1}$ of invert  $\phi^{-1}(v)=(m v-1)/(m-1)$. Fix now $\lambda>0$, with $\bar{u}+\lambda<1$. Expression \eqref{equ_interm_inequ3} may be rewritten as
\begin{align*}
&\Pow(\SU(\W))-  (1-\alpha \pi_0)\bar{u}\\[2pt]
&\leq \E \bigg[  \pi_1 m^{-1}\sum_{i=1}^m F_i\circ \Delta_i(\phi( \wh{u}'_{-i})) -  (1-\alpha \pi_0)\bar{u} \bigg]\\[2pt]
&\leq \pi_1  \P(\Omega_{1}^c)  + G_\W(\bar{u}+\lambda)-\alpha \pi_0 (\bar{u}+\lambda) - (1-\alpha \pi_0)\bar{u} \\[2pt]
&=  \pi_1  \P(\Omega_{1}^c)  + G_\W(\bar{u}+\lambda)- (\bar{u}+\lambda) +  \lambda(1- \alpha \pi_0),
 \end{align*}
 where $\Omega_1$ denotes the event $\big\{\forall i, 1\leq i\leq m, \phi(\wh{u}'_{-i})  \leq \bar{u}+ \lambda \big\}$ and where the last inequality comes from the definition of $\bar{u}$. We upper-bound now the probability of $\Omega_{1}^c$:
 \begin{align}
\prob{\Omega_{1}^c} &\leq \sum_{i=1}^m  \prob{\wh{u}'_{-i}  > \phi^{-1}(\bar{u} + \lambda)}\nonumber\\[2pt]
  &\leq  \sum_{i=1}^m \sum_{u>  \phi^{-1}(\bar{u} + \lambda)} \ind{ u(m-1) \in \{0,1,\ldots,m-1\}}
  \prob{\wh{\mathbb{G}}'_{-i}(u)\geq u}\nonumber\\[2pt]
  &=  \sum_{i=1}^m \sum_{v>  \bar{u} + \lambda} \ind{ vm \in \{1,2,\ldots,m\}} \prob{\wh{\mathbb{G}}'_{-i}(\phi^{-1}(v))\geq \phi^{-1}(v)}\nonumber\\[2pt]
  &\leq  \sum_{i=1}^m \sum_{v>  \bar{u} + \lambda} \ind{ vm \in \{1,2,\ldots,m\}} \prob{\G_\W(v) \geq v - m^{-1}}\nonumber,
  \end{align}
where the last inequality uses that $m \G_\W(v) \geq  (m-1) \G'_{-i}(\phi^{-1}(v))$. As a consequence:
 \begin{align}
\prob{\Omega_{1}^c}
&\leq  m \sum_{\substack{v>  \bar{u} + \lambda \\ vm \in \{1,2,\ldots,m\}}}  \prob{\wh{\mathbb{G}}_{\W}(v)-{{G}}_{\W}(v)\geq v - {{G}}_{\W}(v)-m^{-1}}\nonumber\\
&\leq  m \sum_{\substack{v>  \bar{u} + \lambda \\ vm \in \{1,2,\ldots,m\}}} \prob{\wh{\mathbb{G}}_{\W}(v)-{{G}}_{\W}(v) >( \bar{u}+ \lambda) - {G}_{\W}(\bar{u} + \lambda) - m^{-1}}\nonumber\\
&\leq m^2\exp\big\{-2m\big(( \bar{u}+ \lambda) - {G}_{\W}(\bar{u} + \lambda)- m^{-1} \big)_+^2 \big\},\label{equ2}
\end{align}
where we used successively the assumption in (i) of Proposition~\ref{prop_approxpower} and
Hoeffding's inequality (see \cite{Hoeff1963}) for the last inequality. This finally yields~\eqref{approxpower1}.

The point (ii) of Proposition~\ref{prop_approxpower} is similar: noticing that \eqref{equ_interm_inequ3} is an equality when $\Delta\leq 1$, we obtain
\begin{align*}
&\Pow(\SU(\W))- (1-\alpha \pi_0)\bar{u} \\[2pt]
&\geq - \pi_1  \P(\Omega_{2}^c)   + (G_\W(\bar{u}-\lambda)- (\bar{u}-\lambda ))  -   \lambda(1- \alpha \pi_0),
 \end{align*}
with $\Omega_{2} = \big\{\forall i, 1\leq i\leq m, \phi(\wh{u}'_{-i})  \geq \bar{u} - \lambda \big\}$. Next, we have
\begin{align*}
\P(\Omega_{2}^c) &\leq \sum_{i=1}^m \prob{\wh{u}'_{-i}  < \phi^{-1}(\bar{u} - \lambda)} \\[2pt]
&= \sum_{i=1}^m  \prob{\wh{\mathbb{G}}'_{-i}( \phi^{-1}(\bar{u} - \lambda))< \phi^{-1}(\bar{u} - \lambda)}\\[2pt]
&\leq \sum_{i=1}^m  \prob{\G_\W(\bar{u} - \lambda)< \bar{u} - \lambda},
 \end{align*}
where the last inequality uses that $m \G_\W(\bar{u} - \lambda) \leq  (m-1) \G'_{-i}(\phi^{-1}(\bar{u} - \lambda)) + 1$ and thus $\phi^{-1}(\G_\W(\bar{u} - \lambda))\leq  \G'_{-i}(\phi^{-1}(\bar{u} - \lambda))$. As a consequence, we obtain
\begin{align}
\P(\Omega_{2}^c)&\leq  m\prob{\G_\W(\bar{u} - \lambda)< \bar{u} - \lambda}\nonumber\\[2pt]
&=m  \prob{\wh{\mathbb{G}}_{\W}(\bar{u} - \lambda)-{G}_{\W}(\bar{u} - \lambda )<\bar{u} - \lambda -{G}_{\W}(\bar{u} - \lambda )}\nonumber\\[-1pt]
&\leq m\exp\{-2m(\bar{u} - \lambda -{G}_{\W}(\bar{u} - \lambda ) )_+^2 \},\nonumber
\end{align}
which implies \eqref{approxpower2}.

\subsection[Proof of Theorem~4.3]{Proof of Theorem~\ref{th_asymp}}\label{proof_th_asymp}

First remark that for any weight function sequence of $\mathcal{W}$,
the convergence of $(G_{\W^{(m)}})_m$ to $G^{\infty}$ is uniform,
because all these fonctions are non-decreasing and because
$G^{\infty}$ is assumed to be continuous on $[0,1]$. Next we prove that
$\fix(G_{\W^{(m)}}) \rightarrow \fix(G^{\infty})$. (This will imply
directly that $\fix_\lambda^{+}(G_{\W^{(m)}}) \rightarrow
\fix_\lambda^{+}(G^{\infty})$ and $\fix_\lambda^{-}(G_{\W^{(m)}})
\rightarrow \fix_\lambda^{-}(G^{\infty})$ for
$\lambda<\fix(G^{\infty})$.) For this, take a subsequence $m'$ such
that  $\fix(G_{\W^{(m')}})$ converges and prove that its limit $\ell$
is equal to $\fix(G^{\infty})$. From the uniform convergence and the
continuity of $G^{\infty}$, $\ell$ satisfies $G^{\infty}(\ell)=\ell$.
If $\fix(G^{\infty})=0$, the only possible fixed point of $G^{\infty}$
is $0$ and $\ell=0$. If $\fix(G^{\infty})>0$, $0$ and
$\fix(G^{\infty})$ are the only possible fixed points of $G^{\infty}$
(because $\fix_\lambda^-(G^{\infty})>0$ for
$\lambda<\fix(G^{\infty})$). Next, we have
$G_{\W^{(m')}}(\fix(G^{\infty})/2)\geq \fix(G^{\infty})/2$ for large
$m'$ (because $G^{\infty}(\fix(G^{\infty})/2)> \fix(G^{\infty})/2$ )
and thus $\ell\geq \fix(G^{\infty})/2>0$, which in turn implies
$\ell=\fix(G^{\infty})$.

Fix a sequence of weight vector $(\w^{(m)})_m$ belonging to $\mathcal{W}$. We aim now to prove:
\begin{eqnarray}
\lim_m\Pow(\SU(\W^{\star,(m)})) &=& (1-\pi_0 \alpha)  \lim_m \left\{ \fix(G_{\W^{\star,(m)}})\right\}\label{proof_asymp_optimal_1}\\[2pt]
\lim_m\Pow(\SU(\w^{(m)})) &=& (1-\pi_0 \alpha)  \lim_m \left\{ \fix(G_{\w^{(m)}}) \right\}\label{proof_asymp_optimal_2}
\end{eqnarray}
Expression \eqref{asymp_optimal} will then directly follow from $\fix(G_{\W^{\star,(m)}})\geq \fix(G_{\w^{(m)}})$ (as stated in the proof of Theorem~\ref{Th-Pow-oracle}).

Let us state now \eqref{proof_asymp_optimal_1} (the proof for \eqref{proof_asymp_optimal_2} is similar).  Fix $\lambda>0$ with $\lambda<\pi_0(1-\alpha)$. Applying Proposition~\ref{prop_approxpower}, we obtain that
\begin{align*}
&| \Pow(\SU(\W^{\star,(m)})) - (1-\pi_0 \alpha)  \fix(G_{\W^{\star,(m)}}) | \\
&\leq \varepsilon (m,\fix_\lambda^+({G}_{\W^{\star,(m)}}) ) +  \varepsilon (m,\fix_\lambda^-({G}_{\W^{\star,(m)}}) ) \ind{\lambda < \fix(G_{\W^{\star,(m)}})} - 2\lambda (1-\pi_0 \alpha).
\end{align*}
First, denoting $G^{\star,\infty}=\lim_m G_{\W^{\star,(m)}}$, we have $\fix^+_\lambda(G_{\W^{\star,(m)}}) \rightarrow \fix^+_\lambda(G^{\star, \infty})>0$ and thus $\lim_m \varepsilon (m,\fix_\lambda^+({G}_{\W^{\star,(m)}}) )=0$. Second, if $ \fix(G^{\star,\infty})>0$, we have  $\fix_\lambda^-(G^{\star,\infty})>0$ for $\lambda<\fix(G^{\star,\infty})$ and thus $\lim_m \varepsilon (m,\fix_\lambda^-({G}_{\W^{\star,(m)}}) )= 0$ for $\lambda<\fix(G^{\star,\infty})$. If $ \fix(G^{\star,\infty})=0$, we trivially have that $\ind{\lambda<\fix(G_{\W^{\star,(m)}})}$ is equal to zero for $m$ large. As a result, we obtain  for $\lambda$ small enough that
$$
\limsup_m \left\{ \left| \Pow(\SU(\W^{\star,(m)})) - (1-\pi_0 \alpha)  \fix(G_{\W^{\star,(m)}}) \right| \right\} \leq - 2\lambda (1-\pi_0 \alpha).
$$
This yields \eqref{proof_asymp_optimal_1} by letting $\lambda\rightarrow 0$ and by  noticing that $ \lim_m \left\{ \fix(G_{\W^{\star,(m)}})\right\}$ exists. Finally, we have to check that the use of $\W^{\star,(m)}$ in Proposition~\ref{prop_approxpower} was allowed, i.e. that for all $m$ and $u'>u>u^{\star}:= \fix(G_{\W^{\star,(m)}})$, inequality  $u'-G_{\W^{\star,(m)}}(u')>u-G_{\W^{\star,(m)}}(u)$ holds. For this, we let $\w:=\W^{\star, (m)}(u')$ and $u_\w:=\fix(G_\w)$. Since $u^{\star}\geq u_\w$ and  $G_\w(u_\w)=u_\w$ we have $u'-G_\w(u')>u - G_\w(u)$ ($G_\w$ being strictly concave). Therefore, for this particular weight vector $\w$,
$$
u' -G_{\W^{\star,(m)}}(u') = u' - G_\w(u') > u - G_\w(u) \geq u- G_{\W^{\star,(m)}}(u),
$$
where the last inequality holds because $\W^{\star,(m)}$ is a maximum.

Finally, to get the FDR statement \eqref{asymp_control}, we use the same reasoning as above combined with the following finite FDR approximation result:

\begin{proposition}\label{prop_approxFDR}
In the unconditional model, consider a weight function $\W$ with its associated threshold collection $\Delta$ and put $\bar{u}=\fix(G_\W)$. Assume that for all $u' > u > \bar{u}$, $u'-G_\W(u')>  {u} - G_\W( {u})$ and take $\lambda>0$ with $\lambda<1-\bar{u}$. Then the following bounds
hold:\vadjust{\eject}
\begin{eqnarray}
\FDR(\SU(\W)) &\leq& \pi_0\alpha + \pi_0\alpha m^3 \exp\big\{-2m\big(\fix_\lambda^+(G_\W)- m^{-1} \big)_+^2 \big\}  \nonumber\\
&&+ \pi_0\alpha \ind{\bar{u}>\lambda} \left[ m^2\exp\big\{-2m(\fix_\lambda^-({G}_{\W}) )_+^2 \big\} +\frac{2\lambda}{\bar{u}-\lambda} \right] \nonumber\\
&& +  \pi_0\alpha \ind{\bar{u}\leq \lambda} \left[ m^{-1} \sum_{i=1}^m  \sup_{0<u\leq 2\lambda} \left\{W_i(u)\right\} - 1 \right]. \label{approxFDR1}
\end{eqnarray}
Assuming additionally $\Delta\leq 1$, we have
\begin{eqnarray}
\FDR(\SU(\W)) &\geq& \pi_0\alpha - \pi_0\alpha m^3 \exp\big\{-2m\big(\fix_\lambda^+(G_\W)- m^{-1} \big)_+^2 \big\}  \nonumber\\
&&- \pi_0\alpha \ind{\bar{u}>\lambda} \left[ m^2\exp\big\{-2m(\fix_\lambda^-({G}_{\W}) )_+^2 \big\} +\frac{2\lambda}{\bar{u}+\lambda} \right] \nonumber\\
&&- \pi_0\alpha \ind{\bar{u}\leq \lambda} \left[ 1- m^{-1} \sum_{i=1}^m  \inf_{0<u\leq 2\lambda} \left\{W_i(u)\right\}  \right].  \label{approxFDR2}
\end{eqnarray}
\end{proposition}

To prove Proposition~\ref{prop_approxFDR}, we write the FDR as (using the same reasoning and notation as in Section~\ref{proof_TH3}),
\begin{align}
\label{inequ_approx_FDR}
\FDR(\SU(\W)) &=\pi_0 \sum_{i=1}^m   \E \left[ \frac{\ind{p_i\leq \Delta_i(|R|/m)}}{|R|} \bigg| H_i=0\right]\nonumber\\
&=   \pi_0 \sum_{i=1}^m   \E \left[  \frac{\ind{p_i\leq \Delta_i(\phi(\wh{u}'_{-i}))}}{m\phi(\wh{u}'_{-i})}\bigg| H_i=0\right] \nonumber\\
& \leq  \pi_0 m^{-1} \sum_{i=1}^m    \E \left[  \frac{\Delta_i(\phi(\wh{u}'_{-i}))}{\phi(\wh{u}'_{-i})} \right]\nonumber\\
&\leq  \pi_0 m^{-1} \sum_{i=1}^m   \E \left[  \frac{\Delta_i(\phi(\wh{u}'_{-i}))}{\phi(\wh{u}'_{-i})}
\ind{ \Omega_1 \cap \Omega_2}\right]\nonumber\\
&{}\quad  +  \pi_0 \alpha m 
 (\prob{ \Omega_1^c} + \prob{\Omega_2^c}).
\end{align}
On one hand, when $\bar{u}>\lambda$, we may write
\begin{align}
m^{-1}\sum_{i=1}^m   \E \left[  \frac{\Delta_i(\phi(\wh{u}'_{-i}))}{\phi(\wh{u}'_{-i})} \ind{ \Omega_1 \cap \Omega_2}\right]
&\leq m^{-1} \sum_{i=1}^m  \left[  \frac{\Delta_i(\bar{u}+\lambda)}{\bar{u}-\lambda}\right]  = \alpha + 2\alpha \frac{\lambda}{\bar{u}-\lambda}  \nonumber .
\end{align}
On the other hand, when $\bar{u}\leq \lambda$,  we have $\Omega_2^c=\emptyset$ and
\begin{align}
m^{-1}\sum_{i=1}^m   \E \left[  \frac{\Delta_i(\phi(\wh{u}'_{-i}))}{\phi(\wh{u}'_{-i})} \ind{ \Omega_1}\right]
&\leq \alpha m^{-1} \sum_{i=1}^m  \sup_{0<u\leq2 \lambda} \left\{W_i(u)\right\}  \nonumber .
\end{align}
This implies \eqref{approxFDR1}. The proof for \eqref{approxFDR2} is similar, by noticing that \eqref{inequ_approx_FDR} is an equality when $\Delta\leq 1$.

\subsection[Proof of Proposition~3.2]{Proof of Proposition~\ref{prop:maxpow}}\label{sec:proofweightchoice}

Let $t=\alpha u$ and assume \eqref{assump_class1}--\eqref{assump_class2}--\eqref{assump_class3}. Consider first the conditional model. Following the constrained Lagrange multiplier method, the problem is to maximize in $(\lambda,\w)$ the function:
\begin{equation}
 \nonumber 
L(\lambda,\w)=\sum_{i\in \cH_1} F_i(w_i t)-\lambda\bigg(\sum_{i\in\cH_1}w_i-m\bigg).
\end{equation}
Assume that $(w_i)_i$ is a critical point, i.e. that for each $i$:
$\frac{\partial L}{\partial w_i}(\lambda,\w)=t f_i(w_i t)-\lambda=0,$
so that $w_i= t^{-1} f_i^{-1}(\lambda t^{-1})$. Then, $\lambda$ is chosen such that $\sum_i w_i=m$ i.e. $\lambda=t \Psi^{-1}(t)$ for $\Psi(y)=m^{-1} \sum_{j\in \cH_1} f_j^{-1}(y)$. Hence, we find that the only possible critical point is $\forall i, w_i=W_i^{\star}(u)$. To conclude, it is sufficient to prove that $(W_i^{\star}(u))_i$ is a maximum. The latter holds because for each $i$, $f_i$ is decreasing, so that
$\frac{\partial^2 L}{\partial w_i^2}(\lambda,\w)=t^{2} f'_i(w_i t)<0.$
Therefore, since $\alpha u W_i^{\star}(u) = \zeta^{-1}_i (\alpha u)$, where
$\zeta_i(x)=\Psi( f_i(x)) = m^{-1} \sum_{j\in \cH_1} f_j^{-1}(f_i(x))$
is a differentiable increasing function from $(0,1)$ to $(0,\pi_1)$, we easily check that $\W^{\star}$ satisfies \eqref{equ_hypoweights} and is continuous.  Next, assuming in addition \eqref{assump_class4}, we obtain for all $i$,
\begin{align*}
\lim_{u\rightarrow 0^+} W_i^{\star}(u) &=  \lim_{u\rightarrow 0^+} \frac{\zeta_i^{-1}(\alpha u)}{\alpha u}  = \lim_{\varepsilon \rightarrow 0^+} \frac{\varepsilon}{\zeta_i(\varepsilon)} = \lim_{y \rightarrow f_i(0^+)} \frac{f_i^{-1}(y)}{m^{-1} \sum_{j\in \cH_1} f_j^{-1}(y)},
\end{align*}
which exists in $[0,m]$.
Finally, the results in the unconditional model follow from the same reasoning as above by replacing $\cH_1$ by $\{1,\ldots,m\}$.

\subsection[Proof for the continuous Gaussian case of Section~4.3]{Proof for the continuous Gaussian case of Section~\ref{sec:consis}}\label{sec:proofweightconv}

Fix $u\in(0,1]$. We aim to prove that
\begin{equation}\label{equ:equaim}
\frac{1}{m} \sum_{i=1}^m \ol{\Phi}\bigg(-\frac{\mu(i/m)}{2}+\frac{c^{(m)}(u)}{\mu(i/m)}\bigg)
\rightarrow \int_0^1 \ol{\Phi}\bigg(-\frac{\mu(t)}{2}+\frac{c^{\infty}(u)}{\mu(t)}\bigg) dt .
\end{equation}
First, we have $c^{(m)}(u) \rightarrow c^{\infty}(u)$, because $c^{(m)}(u) =\Psi_m^{-1}(\alpha u)$ and $c^{\infty}(u) =\Psi^{-1}(\alpha u)$ where $\Psi_m(x)= \frac{1}{m} \sum_{i=1}^m \ol{\Phi}\big({\mu(i/m)}/{2}+{x}/{\mu(i/m)}\big)$ and where $\Psi(x)= \int_0^1 \ol{\Phi}\big({\mu(t)}/{2}+{x}/{\mu(t)}\big)dt$ are decreasing continuous functions such that $\Psi_m$ converges uniformly to $\Psi$.
Second, we have for any $\epsilon>0$,
\begin{align*}
&\frac{1}{m} \sum_{i=1}^m \bigg| \ol{\Phi}\bigg(-\frac{\mu(i/m)}{2}+\frac{c^{(m)}(u)}{\mu(i/m)}\bigg)-\ol{\Phi}\bigg(-\frac{\mu(i/m)}{2}+\frac{c^{\infty}(u)}{\mu(i/m)}\bigg)  \bigg| \\
&\leq m^{-1} | \{i\telque \mu(i/m)\leq \epsilon\}| + \frac{ \big| c^{(m)}(u)- c^{\infty}(u) \big| }{\epsilon \sqrt{2\pi}},
\end{align*}
where we used that  $\ol{\Phi}$ is $1/\sqrt{2\pi} $-Lipschitz. Since the measure $m^{-1}\sum_{i=1}^m \delta_{i/m}$ converges weakly to the Lebesgue measure $\Lambda$ on $[0,1]$, we derive the inequality  $\limsup_m (m^{-1} | \{i\telque \mu(i/m)\leq \epsilon\}|) \leq \Lambda(t\telque \mu(t)\leq \epsilon)$. By assumption, $\mu$ is positive on $(0,1]$ so that the latter converges to $0$ as $\epsilon$ tends to $0$. This implies \eqref{equ:equaim}.
\medskip\\

\appendix{}

\noindent{\bf \Large Appendix}

\section{Practical implementation of the new procedures}\label{sec:impl}

\begin{algorithm}(Step-up algorithm for \SU$(\mathbf{W})$) 
\label{algo_MW_stepup}
\begin{itemize}
\item[--] Step $1$: compute for each $i$ the weight vector $(W_i(1))_i$ and the weighted $p$-values $p'_i=p_i/W_i(1)$. If all the weighted $p$-values are less than or equal to $\alpha$, then reject all the null hypotheses. Otherwise go to step $2$.
\item[--] Step $j$ $(j\geq 2)$: put $r=m-j+1$ and $u=r/m$ and compute for each $i$ the weight vector $(W_i(u))_i$ and the weighted $p$-values $p'_i=p_i/W_i(u)$. Order the weighted $p$-values following
$p'_{(1)}\leq \cdots\leq p'_{(m)}$. If $p'_{(r)}\leq \alpha u$, then reject the $r$ null hypotheses corresponding to the smaller weighted $p$-values $p'_{(i)}$, $1\leq i\leq r$. Otherwise go to step $j+1$ (if $j=m$ stop and reject no null hypothesis).
\end{itemize}
\end{algorithm}

\begin{algorithm}(Step-down algorithm for \SD$(\mathbf{W})$) 
\label{algo_MW_stepdown}
\begin{itemize}
\item[--] Step $1$: compute for each $i$ the weight vector $(W_i(1/m))_i$ and the weighted $p$-values $p'_i=p_i/W_i(1/m)$. If the smallest weighted $p$-values is strictly larger than $\alpha/m$, then reject no null hypothesis. Otherwise go to step $2$.
\item[--] Step $j$ $(j\geq 2)$: put $r=j$, $u=r/m$ and compute for each $i$ the weight vector $(W_i(u))_i$ and the weighted $p$-values $p'_i=p_i/W_i(u)$. Order the weighted $p$-values following
$p'_{(1)}\leq \cdots\leq p'_{(m)}$. If $p'_{(r)}>\alpha u$, then reject the $r-1$ null hypotheses corresponding to the smaller weighted $p$-values $p'_{(i)}$, $1\leq i\leq r-1$. Otherwise go to step $j+1$ (if $j=m$ stop and reject all the null hypotheses).
\end{itemize}
\end{algorithm}

\section{Proofs of technical lemmas}\label{sec:lemma_proof}

\textbf{Proof of Lemma~\ref{lemma_stepupdown}.}
Let us first prove the first point. Denote $\Delta'_j(u)= \Delta_j((1-m^{-1})u+m^{-1})$, and $\phi(u):=(1-m^{-1})u+m^{-1}$ with $\phi^{-1}(u)=(m u-1)/(m-1)$,  so that $\Delta_j(u)=\Delta'_j(\phi^{-1 }(u))$. Since $ m \G_\W(u)= (m-1) \G'_{-i}(\phi^{-1}(u)) + \ind{p_i\leq \Delta_i(u)}$ the following equivalence holds when
$p_i\leq \Delta_i(u)$:
\begin{equation}
\G_{\W}(u)\geq u \:\: \Longleftrightarrow\:\:  \wh{\mathbb{G}}'_{-i}(\phi^{-1}(u))\geq \phi^{-1}(u).\label{equ_prooflemma}
\end{equation}
First, assuming $p_i \leq \Delta_i(\wh{u})$, equivalence \eqref{equ_prooflemma} used with $u=\wh{u}$ leads to inequality $ \wh{\mathbb{G}}'_{-i}(\phi^{-1}(\wh{u}))\geq \phi^{-1}(\wh{u})$ and thus  $\phi^{-1}(\wh{u})\leq \wh{u}'_{-i}$ because $\wh{u}'_{-i}$ is defined as a maximum. This implies $p_i \leq  \Delta_i(\wh{u})= \Delta'_i(\phi^{-1}(\wh{u})) \leq  \Delta'_i(\wh{u}'_{-i})$.
Conversely, assuming  $p_i \leq \Delta'_i(\wh{u}'_{-i}) = \Delta_i( \phi(\wh{u}'_{-i}))$, equivalence \eqref{equ_prooflemma} used with $u= \phi(\wh{u}'_{-i})$ yields $\G_\W( \phi(\wh{u}'_{-i}))\geq  \phi(\wh{u}'_{-i})$ and thus $ \phi(\wh{u}'_{-i}) \leq \wh{u}$ by definition of $\wh{u}$. The first point is thus proved by additionally noticing that when both $p_i \leq \Delta_i(\wh{u})$ and $p_i \leq \Delta'_i(\wh{u}'_{-i})$ we have both $\wh{u}'_{-i}\geq \phi^{-1}(\wh{u})$  and $ \phi(\wh{u}'_{-i}) \leq \wh{u}$, so that  $\phi(\wh{u}'_{-i})=  \wh{u}$. 

For the second point of the lemma, remark that
\begin{equation}
m \G_{\W}(u)= (m-1) \G_{-i}( um/(m-1) ) + \ind{p_i\leq \Delta_i(u)}.\label{equ_lemma_useful}
\end{equation}
Therefore, we always have
$
\G_{\W}(u) \geq u \:\: \Longleftarrow\:\:   \G_{-i} ( um/(m-1) )\geq  um/(m-1),
$
which implies, using  $u=(1-m^{-1})\wh{u}_{-i} $, that $\wh{u} \geq (1-m^{-1})\wh{u}_{-i} $ always holds. Next,  when $p_i > \Delta_i(u)$, we have
$
\G_{\W}(u) \geq u \:\: \Longrightarrow\:\:   \G_{-i} ( um/(m-1) )\geq  um/(m-1),
 $
so that taking $u=\wh{u}$ in the relation above  leads to $\wh{u} \leq (1-m^{-1})\wh{u}_{-i} $ and thus $\wh{u} = (1-m^{-1})\wh{u}_{-i} $. Conversely, if $p_i \leq \Delta_i(\wh{u})$, from the first point of the lemma we obtain $\wh{u} = (1-m^{-1})\wh{u}'_{-i}+m^{-1} $ and since $\wh{u}'_{-i}\geq \wh{u}_{-i}$ (because pointwise $\G_{-i}\leq \G'_{-i}$), we deduce  $\wh{u} > (1-m^{-1})\wh{u}_{-i}$ which finishes the proof.
\qed\\

\noindent\textbf{Proof of Lemma~\ref{lemma_stepdown}.}
For the first point, take $u'\leq  (k-1)/(m-1)$ and apply \eqref{equ_lemma_useful} with $u=(1-m^{-1})u'$, which gives  $m \G_{\W}((1-m^{-1})u')= (m-1) \G_{-i}( u' ) + \ind{p_i\leq \Delta_i((1-m^{-1})u')}$. Since $(1-m^{-1})u'\leq  (k-1)/m$, assuming $\wt{u}\geq k/m$ and $p_i>\Delta_i((k-1)/m)$, we obtain $p_i> \Delta_i((1-m^{-1})u')$ and $ \G_{\W}((1-m^{-1})u')\geq (1-m^{-1})u'$, which thus leads to $ \G_{-i}( u' )\geq u'$. Since this holds for any $u'\leq  (k-1)/(m-1)$, we finally have $\wt{u}_{-i} \geq (k-1)/(m-1)$.

To prove the second point, take $u'\leq (k-1)/m$ and use \eqref{equ_lemma_useful} with $u=u'$. This gives $m \G_{\W}(u')= (m-1) \G_{-i}( u'm/(m-1) ) + \ind{p_i\leq \Delta_i(u')}$. Since $u'm/(m-1)\leq (k-1)/(m-1)$ and $u'\leq (k-1)/m$, if $\wt{u}_{-i} \geq (k-1)/(m-1)$ and $p_i\leq \Delta_i((k-1)/m)$, we obtain  $\G_{\W}(u') \geq u' + m^{-1}\geq u'$. This holds for any $u'\leq (k-1)/m$ and also for $u'=k/m$ because $\G_{\W}(k/m)\geq \G_{\W}((k-1)/m) \geq (k-1)/m + m^{-1} = k/m$. Finally $\wt{u}\geq k/m$.

For the third point, remark that since we trivially have $ \wt{u} \geq (1-m^{-1})\wt{u}_{-i}$ (with an argument similar than in the step-up case), it is sufficient to prove $ \wt{u} \leq (1-m^{-1})\wt{u}_{-i}$. For this, use \eqref{equ_lemma_useful} with $u'=(1-m^{-1})\wt{u}_{-i} + m^{-1}$, leading to $m \G_{\W}(u')= (m-1) \G_{-i}( u'm/(m-1) ) + \ind{p_i\leq \Delta_i(u')}$. Since $R_{-i}$ is step-down and since $u'm/(m-1)= \wt{u}_{-i} + (m-1)^{-1} $ we have by definition of $ \wt{u}_{-i} $ that $\G_{-i}(u'm/(m-1))> u'm/(m-1)$. Therefore, assuming $p_i>\Delta_i(u')$, we obtain that $ \G_{\W}(u')>  u'$, meaning that $ \wt{u} < u'$ because $R$ is step-down. Hence $ \wt{u} \leq   (1-m^{-1})\wt{u}_{-i}$.
\qed 

\section{Some FDR bounds for \SU$(\mathbf{W})$ and \SD$(\mathbf{W})$}\label{appen_contrexFDRcontrol}

\subsection{Step-up case}

\begin{lemma}
Consider the conditional model in the situation where only two true hypotheses are tested, that is, $m=m_0=2$. Then for any weight function $W$, the procedure \SU$(\mathbf{W})$ has a FDR equal to $\alpha+\alpha^2 (1-W_1(1))(W_1(1)-W_1(1/2)).$
\end{lemma}

In particular, in the conditional model with $m=m_0=2$, the FDR in the above lemma has a maximum equal to $\alpha+\alpha^2/4$, attained e.g. for the weight function $W_1(1/2)=0$; $W_2(1/2)=2$; $W_1(1)=0.5$; $W_2(1)=1.5$.
Additionally, in the unconditional model, the FDR of the above procedure is larger than $\pi_0^2 (\alpha+\alpha^2/4)$ and can thus be larger than $\alpha$ when $\pi_0$ is (very) close to $1$.
\subsection{Step-down case}

The next result states that $\SD(\W)$ control non-asymptotically the FDR without correction in the case $m=2$ and when $m_0=m$ in the conditional model. This is quite intriguing and we may think that $\SD(\W)$ controls the FDR for any $m$ and $m_0$.

\begin{lemma}
For any weight function $\W$, the procedure \SD$(\mathbf{W})$ controls the FDR at level $\alpha$ in either of the two following cases:
\begin{itemize}
\item[(i)] in the unconditional model when all the hypotheses are true, that is $m_0=m$,
\item[(ii)] in both conditional and unconditional model when $m=2$. 
\end{itemize}
\end{lemma}

To prove (i), we easily check that, when all the hypotheses are true,
the FDR of \SD$(\mathbf{W})$ is $1-\probnn{\G_{W}(1/m)=0}$ and is thus equal to the
FDR of {\bf LSD}$(\mathbf{W}(1/m))$, which is equal to $\alpha$ from results on weighted linear step down procedures (see e.g. \cite{BR2008EJS}). To prove point (ii), we just have to check the case $m_0=1$ from point (i). This trivially holds from \eqref{equ_boundFDRMWLSU} (which also holds in the step-down case), because all the weights are smaller than $m=2$.


\begin{thebibliography}{}

\bibitem[Benjamini and Hochberg, 1995]{BH1995}
\textsc{Benjamini, Y. and Hochberg, Y.} (1995).
\newblock Controlling the false discovery rate: a practical and powerful
  approach to multiple testing.
\newblock \textit{J. Roy. Statist. Soc. Ser. B}, 57(1):289--300.
\MR{1325392}

\bibitem[Benjamini et~al., 2006]{BKY2006}
\textsc{Benjamini, Y., Krieger, A.M., and Yekutieli, D.} (2006).
\newblock Adaptive linear step-up procedures that control the false discovery
  rate.
\newblock \textit{Biometrika}, 93(3):491--507.
\MR{2261438}

\bibitem[Benjamini and Yekutieli, 2001]{BY2001}
\textsc{Benjamini, Y. and Yekutieli, D.} (2001).
\newblock The control of the false discovery rate in multiple testing under
  dependency.
\newblock \textit{Ann. Statist.}, 29(4):1165--1188.
\MR{1869245}

\bibitem[Blanchard and Roquain, 2008]{BR2008EJS}
\textsc{Blanchard, G. and Roquain, E.} (2008).
\newblock Two simple sufficient conditions for {FDR} control.
\newblock \textit{Electron. J. Stat.}, 2:963--992.
\MR{2448601}

\bibitem[Blanchard and Roquain, 2009]{BR2008b}
\textsc{Blanchard, G. and Roquain, E.} (2009).
\newblock Adaptive false discovery rate control under independence and dependence.
\newblock \textit{J. Mach. Learn. Res.}
\newblock To appear.

\bibitem[Efron et~al., 2001]{ETST2001}
\textsc{Efron, B., Tibshirani, R., Storey, J.D., and Tusher, V.}
(2001).
\newblock Empirical {B}ayes analysis of a microarray experiment.
\newblock \textit{J. Amer. Statist. Assoc.}, 96(456):1151--1160.
\MR{1946571}

\bibitem[Finner et~al., 2009]{FDR2009}
\textsc{Finner, H., Dickhaus, R., and Roters, M.} (2009).
\newblock On the false discovery rate and an asymptotically optimal rejection
  curve.
\newblock \textit{Ann. Statist.}, 37(2):596--618.
\MR{2502644}

\bibitem[Gavrilov et~al., 2009]{GBS2008}
\textsc{Gavrilov, Y., Benjamini, Y., and Sarkar, S.K.} (2009).
\newblock An adaptive step-down procedure with proven fdr control under
  independence.
\newblock \textit{Ann. Statist.}, 37(2):619--629.
\MR{2502645}

\bibitem[Genovese and Wasserman, 2004]{GW2004}
\textsc{Genovese, C. and Wasserman, L.} (2004).
\newblock A stochastic process approach to false discovery control.
\newblock \textit{Ann. Statist.}, 32(3):1035--1061.
\MR{2065197}

\bibitem[Genovese et~al., 2006]{GRW2006}
\textsc{Genovese, C.R., Roeder, K., and Wasserman, L.} (2006).
\newblock False discovery control with {$p$}-value weighting.
\newblock \textit{Biometrika}, 93(3):509--524.
\MR{2261439}

\bibitem[Hoeffding, 1963]{Hoeff1963}
\textsc{Hoeffding, W.} (1963).
\newblock Probability inequalities for sums of bounded random variables.
\newblock \textit{J. Amer. Statist. Assoc.}, 58:13--30.
\MR{0144363}

\bibitem[Holm, 1979]{Holm1979}
\textsc{Holm, S.} (1979).
\newblock A simple sequentially rejective multiple test procedure.
\newblock \textit{Scand. J. Statist.}, 6(2):65--70.
\MR{0538597}

\bibitem[Hyman et~al., 2002]{Hyman2002}
\textsc{Hyman, E., Kauraniemi, P., Hautaniemi, S., Wolf, M., Mousses,
S., Rozenblum,
  E., Ringnr, M., Sauter, G., Monni, O., Elkahloun, A., Kallioniemi, A., and
  Kallioniemi, O.} (2002).
\newblock Impact of dna amplification on gene expression patterns in breast
  cancer.
\newblock \textit{Cancer research}, 62(21):6240--5.

\bibitem[Kim and van~de Wiel, 2008]{KW2008}
\textsc{Kim, K.I. and van~de Wiel, M.} (2008).
\newblock Effects of dependence in high-dimensional multiple testing problems.
\newblock \textit{BMC Bioinformatics}, 9(1):114.

\bibitem[Lehmann et~al., 2005]{LRS2005}
\textsc{Lehmann, E.L., Romano, J.P., and Shaffer, J.P.} (2005).
\newblock On optimality of stepdown and stepup multiple test procedures.
\newblock \textit{Ann. Statist.}, 33(3):1084--1108.
\MR{2195629}

\bibitem[Muris et~al., 2007]{Muris2007}
\textsc{Muris, J., Ylstra, B., Cillessen, S., Ossenkoppele, G.,
Kluin-Nelemans, J.,
  Eijk, P., Nota, B., Tijssen, M., de~Boer, W., van~de Wiel, M., van~den
  Ijssel, P., Jansen, P., de~Bruin, P., van Krieken,~J., Meijer, G., Meijer,
  C., and Oudejans, J.} (2007).
\newblock {{P}rofiling of apoptosis genes allows for clinical stratification of
  primary nodal diffuse large {B}-cell lymphomas}.
\newblock \textit{Br. J. Haematol.}, 136:38--47.

\bibitem[Picard et~al., 2007]{Picard2007}
\textsc{Picard, F., Robin, S., Lebarbier, E., and Daudin, J.} (2007).
\newblock {{A} segmentation/clustering model for the analysis of array
  {C}{G}{H} data}.
\newblock \textit{Biometrics}, 63:758--766.
\MR{2395713}

\bibitem[Romano et~al., 2008]{RSW2008}
\textsc{Romano, J.P., Shaikh, A.M., and Wolf, M.} (2008).
\newblock Control of the false discovery rate under dependence using the
  bootstrap and subsampling.
\newblock \textit{TEST: An Official Journal of the Spanish Society of Statistics
  and Operations Research}, 17(3):417--442.
\MR{2470085}

\bibitem[Roquain, 2007]{Roq2007}
\textsc{Roquain, E.} (2007).
\newblock \textit{Exceptional motifs in heterogeneous sequences. Contributions to
  theory and methodology of multiple testing}.
\newblock PhD thesis, Universit\'e Paris XI.

\bibitem[Roquain and van~de Wiel, 2008]{RW2008}
\textsc{Roquain, E. and van~de Wiel, M.} (2008).
\newblock Multi-weighting for {FDR} control.
\newblock {ArXiV} preprint math.ST/0807.4081v1.

\bibitem[Rubin et~al., 2006]{RDV2006}
\textsc{Rubin, D., Dudoit, S., and van~der Laan, M.} (2006).
\newblock A method to increase the power of multiple testing procedures through
  sample splitting.
\newblock \textit{Stat. Appl. Genet. Mol. Biol.}, 5:Art. 19, 20 pp. (electronic).
\MR{2240850}

\bibitem[Sarkar, 2002]{Sar2002}
\textsc{Sarkar, S.K.} (2002).
\newblock Some results on false discovery rate in stepwise multiple testing
  procedures.
\newblock \textit{Ann. Statist.}, 30(1):239--257.
\MR{1892663}

\bibitem[Storey, 2002]{Storey2002}
\textsc{Storey, J.D.} (2002).
\newblock A direct approach to false discovery rates.
\newblock \textit{J. R. Stat. Soc. Ser. B Stat. Methodol.}, 64(3):479--498.
\MR{1924302}

\bibitem[Storey, 2003]{Storey2003}
\textsc{Storey, J.D.} (2003).
\newblock The positive false discovery rate: a {B}ayesian interpretation and
  the {$q$}-value.
\newblock \textit{Ann. Statist.}, 31(6):2013--2035.
\MR{2036398}

\bibitem[Storey, 2007]{Storey2007}
\textsc{Storey, J.D.} (2007).
\newblock The optimal discovery procedure: a new approach to simultaneous
  significance testing.
\newblock \textit{J. R. Stat. Soc. Ser. B Stat. Methodol.}, 69(3):347--368.
\MR{2323757}

\bibitem[van~der Vaart, 1998]{Vaart1998}
\textsc{van~der Vaart, A.W.} (1998).
\newblock \textit{Asymptotic statistics}, volume~3 of \textit{Cambridge Series in
  Statistical and Probabilistic Mathematics}.
\newblock Cambridge University Press, Cambridge.
\MR{1652247}

\bibitem[Wasserman and Roeder, 2006]{WR2006}
\textsc{Wasserman, L. and Roeder, K.} (2006).
\newblock Weighted hypothesis testing.
\newblock Technical report, Dept. of statistics, Carnegie Mellon University.

\end{thebibliography}
\bibliographystyle{apalike}

\end{document}